\documentclass[12pt]{article}
\usepackage{amssymb}

\usepackage[latin1]{inputenc}
\usepackage{graphics}
\usepackage{graphicx}
\usepackage{epsfig}
\usepackage{amsfonts}
\usepackage{amscd}
\usepackage{color}
\usepackage{shadow}
\usepackage{a4wide}
\usepackage[all]{xy}
\usepackage{latexsym}
\usepackage{graphicx}
\usepackage{amsmath}
\usepackage{hyperref,amsthm}

\numberwithin{equation}{section}

\normalsize \makeatletter

\def\cdotv{\raise 2pt\hbox{,}}

 \def\pasdegrille{\let\grille =
\pasgrille}  \def\aat#1#2#3{ \divide
\dimen1 by 48 \dimen3=\dimen1 \multiply \dimen1 by #1 \advance \dimen1
by -\dimen3 \divide \dimen1 by 101 \multiply \dimen1 by 100 \divide
\dimen2 by \count11 \multiply \dimen2 by #2
\setbox0=\hbox{#3}\ht0=0pt\dp0=0pt \rlap{\kern\dimen1 \vbox
to0pt{\kern-\dimen2\box0\vss}}\dimen1= \wd1 \dimen2=\ht1}
\def\pasgrille{ \count12= \dimen1 \divide \count12 by 50 \divide
\dimen2 by \count12 \count11 =\dimen2 \ \divide \dimen1 by 48
\setlength{\unitlength}{\dimen1} \smash{\rlap{\ }} \dimen1= \wd1
\dimen2=\ht1 } \def\grille{ \count12= \dimen1 \divide \count12 by 50
\divide \dimen2 by \count12 \count11 =\dimen2 \ \divide \dimen1 by 48
\setlength{\unitlength}{\dimen1} \smash{\rlap{\graphpaper[1](0,0)(50,
\count11)}} \dimen1= \wd1 \dimen2=\ht1 }

\edef\@tempa#1#2{\def#1{\mathaccent\string"\noexpand\accentclass@#2 }}
\@tempa\ring{017}

\newcommand{\be}{\begin{equation}}
\newcommand{\ee}{\end{equation}}
\newcommand{\bea}{\begin{eqnarray}}
\newcommand{\eea}{\end{eqnarray}}

\newtheorem{thm}{Theorem}[section]
\newtheorem{cor}[thm]{Corollary}
\newtheorem{lemma}[thm]{Lemma}
\newtheorem{prop}[thm]{Proposition}
\newtheorem{dfn}[thm]{Definition}
\theoremstyle{remark}
\newtheorem{rmq}[thm]{Remark}
\newtheorem{assum}[thm]{Assumptions}

\begin{document}

\title{On the Schr\"{o}dinger equation outside strictly convex obstacles}
\date{}
\author{Oana Ivanovici  \footnote{The author was supported by grant A.N.R.-07-BLAN-0250}
 \\Universite Paris-Sud, Orsay,\\
Mathematiques, Bat. 430, 91405 Orsay Cedex, France\\
oana.ivanovici@math.u-psud.fr}

\maketitle \baselineskip 17pt
\begin{abstract}
We prove sharp Strichartz estimates for the semi-classical Schr\"{o}dinger equation on a compact Riemannian manifold with smooth, strictly geodesically concave boundary. We deduce classical Strichartz estimates for the Schr\"odinger equation outside a strictly convex obstacle, local existence for the $H^{1}$-critical (quintic) Schr\"{o}dinger equation and scattering for the sub-critical Schr\"odinger equation in $3D$.
\end{abstract}

\section{Introduction}
Let $(M,g)$ be a Riemannian manifold of dimension $n\geq 2$. Strichartz estimates are a family of dispersive estimates on solutions $u(x,t):M\times[-T,T]\rightarrow\mathbb{C}$ to the Schr\"{o}dinger equation
\be
i\partial_{t}u+\Delta_{g}u=0,\quad u(x,0)=u_{0}(x),
\ee
where $\Delta_{g}$ denotes the Laplace-Beltrami operator on $(M,g)$. In their most general form, local Strichartz estimates state that 
\be
\|u\|_{L^{q}([-T,T],L^{r}(M))}\leq C\|u_{0}\|_{H^{s}(M)},
\ee
where $H^{s}(M)$ denotes the Sobolev space over $M$ and $2\leq q,r\leq\infty$ satisfy $(q,r,n)\neq (2,\infty,2)$ (for the case $q=2$ see \cite{keta98}) and are given by the scaling admissibility condition
\be\label{scaling}
\frac{2}{q}+\frac{n}{r}=\frac{n}{2}.
\ee
In $\mathbb{R}^{n}$ and for $g_{ij}=\delta_{ij}$, Strichartz estimates in the context of the wave and Schr\"{o}dinger equations have a long history, beginning with Strichartz pioneering work \cite{stri77}, where he proved the particular case $q=r$ for the wave and (classical) Schr\"{o}dinger equations. This was later generalized to mixed $L^{q}_{t}L^{r}_{x}$ norms by Ginibre and Velo \cite{give85} for Schr\"{o}dinger equations, where $(q,r)$ is sharp admissible and $q>2$; the wave estimates were obtained independently by Ginibre-Velo \cite{give95} and Lindblad-Sogge \cite{ls95}, following earlier work by Kapitanski \cite{lev90}. The remaining endpoints for both equations  were finally settled by Keel and Tao \cite{keta98}. In that case $s=0$ and $T=\infty$; (see also Kato \cite{ka87}, Cazenave-Weissler \cite{cawe90}). Estimates for the flat $2$-torus were shown by Bourgain \cite{bourg03} to hold for $q=r=4$ and any $s>0$.

In the variable coefficients case, even without boundaries, the situation is much more complicated:  we simply recall the pioneering work of Staffilani and Tataru \cite{stta02}, dealing with compact, non trapping pertubations of the flat metric, the works by Hassell, Tao and Wunsch \cite{hatawu06}, by Robbiano and Zuily \cite{rozu05} and the recent work of Bouclet and Tzvetkov \cite{botz} which considerably weakens the decay of the pertubation (retaining the non-trapping character at spatial infinity). On compact manifolds without boundaries, Burq, Gerard and Tzvetkov \cite{bgt04} established Strichartz estimates with $s=1/p$, hence with a loss of derivatives when compared to the case of flat geometries. Recently, Blair, Smith and Sogge \cite{bss08} improved on the current results for compact $(M,g)$ where either $\partial M\neq\emptyset$, or $\partial M=\emptyset$ and $g$ Lipschitz, by showing that Strichartz estimates hold with a loss of $s=4/3p$ derivatives. This appears to be the natural analog of the estimates of \cite{bgt04} for the general boundaryless case.
 
In this paper we prove that Strichartz estimates for the semi-classical Schr\"{o}dinger equation also hold on Riemannian manifolds with smooth, strictly geodesically concave boundaries. By the last condition we understand that the second fundamental form on the boundary of the manifold is strictly positive definite. 
We have two main examples of such manifolds in mind: first, we consider the case of a compact manifold with strictly concave boundary, which we shall denote $S$ in the rest of the paper. 
The second example is the exterior of the strictly convex obstacle in $\mathbb{R}^{n}$, which will be denoted by $\Omega$.

\begin{assum}\label{assum1}
Let $(S,g)$ be a smooth $n$-dimensional compact Riemannian manifold with $C^{\infty}$ boundary. Assume that $\partial S$ is strictly geodesically concave. Let $\Delta_{g}$ be the Laplace-Beltrami operator associated to $g$.

Let $0<\alpha_{0}\leq1/2$, $2\leq\beta_{0}$, $\Psi\in C^{\infty}_{0}(\mathbb{R}\setminus\{0\})$ be compactly supported in the interval $(\alpha_{0},\beta_{0})$. We introduce the operator $\Psi(-h^{2}\Delta_{g})$ using the Dynkin-Helffer-Sj\"ostrand formula \cite{Dav89} and refer to \cite{nier93}, \cite{Dav89} or \cite{ivpl08} for a complete overview of its properties (see also \cite{bgt04} for compact manifolds without boundaries).
\begin{dfn}\label{dfnpsi}
Given $\Psi\in C^{\infty}_{0}(\mathbb{R})$ we have
\[
\Psi(-h^{2}\Delta_{g})=-\frac{1}{\pi}\int_{\mathbb{C}}\bar{\partial}\tilde{\Psi}(z)(z+h^{2}\Delta_{g})^{-1}dL(z),
\]
where $dL(z)$ denotes the Lebesque measure on $\mathbb{C}$ and $\tilde{\Psi}$ is an almost analytic extension of $\Psi$, e.g., with $<z>=(1+|z|^{2})^{1/2}$, $N\geq 0$,
\[
\tilde{\Psi}(z)=\Big(\sum_{m=0}^{N}\partial^{m}\Psi(\text{Re}z)(i\text{Im}z)^{m}/m!\Big)\tau(\text{Im}z/<\text{Re}z>),
\]
where $\tau$ is a non-negative $C^{\infty}$ function such that $\tau(s)=1$ if $|s|\leq 1$ and $\tau(s)=0$ if $|s|\geq 2$.
\end{dfn}
\end{assum}

Our main result reads as follows:
\begin{thm}\label{thmstrichartz}
Under the Assumptions \ref{assum1}, given $(q,r)$ satisfying the scaling condition \eqref{scaling}, $q>2$ and $T>0$ sufficiently small, there exists a constant $C=C(T)>0$ such that the solution $v(x,t)$ of the semi-classical Schr\"{o}dinger equation on $S\times\mathbb{R}$ with Dirichlet boundary condition
\be\label{schrodscl}
\left\{\begin{array}{c}
ih\partial_{t}v+h^{2}\Delta_{g}v=0 \quad \text{on} \quad S\times\mathbb{R},\\
v(x,0)=\Psi(-h^{2}\Delta_{g})v_{0}(x),\\
v|_{\partial S}=0\\
\end{array}
\right.
\ee
satisfies
\be\label{sescl}
\|v\|_{L^{q}((-T,T),L^{r}(S))}\leq Ch^{-\frac{1}{q}}\|\Psi(-h^{2}\Delta_{g})v_{0}\|_{L^{2}(S)}.
\ee
\end{thm}
\begin{rmq}
An example of compact manifold with smooth, strictly concave boundary is given by the Sina\"i billiard (defined as the complementary of a strictly convex obstacle on a cube of $\mathbb{R}^{n}$ with periodic boundary conditions). 
\end{rmq}
We deduce from Theorem \ref{thmstrichartz} and \cite[Thm.1.1]{ivpl08} (see also Lemma \ref{lemf}), as in \cite{bgt04}, the following Strichartz estimates with derivative loss:
\begin{cor}
Under the Assumptions \ref{assum1}, given $(q,r)$ satisfying the scaling condition \eqref{scaling}, $q>2$ and $I$ any finite time interval, there exists a constant $C=C(I)>0$ such that the solution $u(x,t)$ of the (classical) Schr\"{o}dinger equation on $S\times\mathbb{R}$ with Dirichlet boundary condition 
\be\label{schrodcl}
\left\{\begin{array}{c}
i\partial_{t}u+\Delta_{g}u=0, \quad  \text{on} \quad S\times\mathbb{R},\\
u(x,0)=u_{0}(x),\quad
u|_{\partial S}=0
\end{array}
\right.
\ee
satisfies
\be
\|u\|_{L^{q}((I,L^{r}(S))}\leq C(I)\|u_{0}\|_{H^{\frac{1}{q}}(S)}.
\ee
\end{cor}

The proof of Theorem \ref{thmstrichartz} is based on the finite speed of propagation of the semi-classical flow (see Lebeau \cite{le92}) and the energy conservation which allow us to use the arguments of Smith and Sogge \cite{smso95} for the wave equation: using the Melrose and Taylor parametrix for the stationary wave (see \cite{meltay85}, \cite{meta85} or Zworski \cite{zw90}) we obtain, by Fourier transform in time, a parametrix for the Schr\"odinger operator near a "glancing" point. Since in the elliptic and hyperbolic regions the solution of \eqref{schrod} will clearly satisfy the same Strichartz estimates as on a manifold without boundary (in which case we refer to \cite{bgt04}), we need to restrict our attention only on the glancing region.

As an application of Theorem \ref{thmstrichartz} we prove classical, global Strichartz estimates for the Schr\"odinger equation outside a strictly convex domain in $\mathbb{R}^{n}$.

\begin{assum}\label{assum2} 
Let $\Omega=\mathbb{R}^{n}\setminus\Theta$, where $\Theta$ is a compact with smooth boundary. We assume that $n\geq 2$ and that $\partial\Omega$ is strictly geodesically concave throughout. Let 
$\Delta_{D}=\sum_{j=1}^{n}\partial^{2}_{j}$
denote the Dirichlet Laplace operator (with constant coefficients) on $\Omega$.
\end{assum}
\begin{thm}\label{thmstri}
Under the Assumptions \ref{assum2}, given $(q,r)$ satisfying the scaling condition \eqref{scaling}, $q>2$ and $u_{0}\in L^{2}(\Omega)$, there exists a constant $C>0$ such that the solution $u(x,t)$ of the Schr\"{o}dinger equation on $\Omega\times\mathbb{R}$ with Dirichlet boundary condition
\be\label{schrod}
\left\{\begin{array}{c}
i\partial_{t}u+\Delta_{D}u=0, \ on \ \Omega\times\mathbb{R},\\
u(x,0)=u_{0}(x),\\
u|_{\partial\Omega}=0\\
\end{array}
\right.
\ee
satisfies
\be\label{se}
\|u\|_{L^{q}(\mathbb{R},L^{r}(\Omega))}\leq C\|u_{0}\|_{L^{2}(\Omega)}.
\ee
\end{thm}
The proof of Theorem \ref{thmstri} combines several arguments: firstly, we perform a time rescaling, first used by Lebeau \cite{le92} in the context of control theory, which transforms the equation into a semi-classical problem for which we can use the local in time semi-classical Strichartz estimates proved in Theorem \ref{thmstrichartz}.
Secondly, we adapt a result of Burq \cite{bu02} which provides Strichartz estimates without loss for a non-trapping problem, with a metric that equals the identity outside a compact set. The proof relies on a local smoothing effect for the free evolution $\exp{(it\Delta_{D})}$, first observed  independently by Constantin and Saut \cite{cosa89}, Sj\"olin \cite{Sj87} and Vega \cite{veg88} in the flat case, and then by Doi \cite{doi96} on non-trapping manifolds and by Burq, G\'erard and Tzvetkov \cite{bgt03} on exterior domains. Following a strategy suggested by Staffilani and Tataru \cite{stta02}, we prove that away from the obstacle the free evolution enjoys the Strichartz estimates exactly as for the free space. 

We give two applications of Theorem \ref{thmstri} : the first one is a local existence result for the quintic Schr\"{o}dinger equation in $3D$, while the second one is a scattering result for the subcritical (sub-quintic) Schr\"odinger equation in $3D$ domains.

\begin{thm}\label{thmgeq}(Local existence for the quintic Schr\"{o}dinger equation)
Let $\Omega$ be a three dimensional Riemannian manifold satisfying the Assumptions \ref{assum2}. Let $T>0$ and $u_{0}\in H^{1}_{0}(\Omega)$. Then there exists a unique solution $u\in C([0,T],H^{1}_{0}(\Omega))\cap L^{5}((0,T],W^{1,30/11}(\Omega))$ of the quintic nonlinear equation
\begin{equation}\label{eqquint}
i\partial_{t}u+\Delta_{D}u=\pm |u|^{4}u\ on\ \Omega\times\mathbb{R},\quad u|_{t=0}=u_{0}\ on \ \Omega, \quad u|_{\partial\Omega}=0.
\end{equation}
Moreover, for any $T>0$, the flow $u_{0}\rightarrow u$ is Lipschitz continuous from any bounded set of $H^{1}_{0}(\Omega)$ to $C([-T,T),H^{1}_{0}(\Omega))$. If the initial data $u_{0}$ has sufficiently small $H^{1}$ norm, then the solution is global in time.
\end{thm}

\begin{thm}\label{thmscattering}(Scattering for subcritical Schr\"{o}dinger equation) Let $\Omega$ be a three dimensional Riemannian manifold satisfying the Assumptions \ref{assum2}. Let $1+\frac{4}{3}\leq p<5$ and $u_{0}\in H^{1}_{0}(\Omega)$. Then the global in time solution of the defocusing Schr\"{o}dinger equation
\be\label{scds}
i\partial_{t}u+\Delta_{D}u=|u|^{p-1}u,\quad u|_{t=0}=u_{0}\ on\ \Omega,\quad u|_{\partial\Omega}=0
\ee 
scatters in $H^{1}_{0}(\Omega)$. If $p=5$ and the gradient $\nabla u_{0}$ of the initial data has sufficiently small $L^{2}$ norm, then the global solution of the critical Schr\"{o}dinger equation scatters in $H^{1}_{0}(\Omega)$.
\end{thm}
Results for the Cauchy problem associated to the critical wave equation outside a strictly convex obstacle were obtained by Smith and Sogge \cite{smso95}. Their result was a consequence of the fact that the Strichartz estimates for the Euclidian wave equation also hold on Riemannian manifolds with smooth, compact and strictly concave boundaries. 

In \cite{bulepl07}, Burq, Lebeau and Planchon proved that the defocusing quintic wave equation with Dirichlet boundary conditions is globally wellposed on $H^{1}(M)\times L^{2}(M)$ for any smooth, compact domain $M\subset\mathbb{R}^{3}$. Their proof relies on $L^{p}$ estimates for the spectral projector obtained by Smith and Sogge \cite{smso06}. A similar result for the defocusing critical wave equation with  Neumann boundary conditions was obtained in \cite{bupl07}. 

In the case of Schr\"odinger equation in $\mathbb{R}^{3}\times\mathbb{R}_{t}$, Colliander, Keel, Staffilani, Takaoka and Tao \cite{ckstt06} established global well-posedness and scattering for energy-class solutions to the quintic defocusing Schr\"{o}dinger equation \eqref{eqquint}, which is energy-critical. When the domain is the complementary of an obstacle in $\mathbb{R}^{3}$, non-trapping but not convex, the counterexamples constructed in \cite{ivond08} for the wave equation suggest that losses are likely to occur in the Strichartz estimates for the Schr\"odinger equation too. In this case
 Burq, Gerard and Tzvetkov \cite{bgt03} proved global existence for subcubic defocusing nonlinearities and Anton \cite{an99} for the cubic case. Recently, Planchon and Vega \cite{plve08} improved the local well-posedness theory to $H^{1}$-subcritical (subquintic) nonlinearities for $n=3$. 
Theorem \ref{thmscattering} is proved in \cite{plve08} in the case of the exterior of a star-shaped domain for the particular case $p=3$, using the following estimate on the solution to the linear problem
\[
\|u\|^{4}_{L^{4}_{t,x}}\lesssim\|u_{0}\|^{3}_{L^{2}}\|\nabla u_{0}\|_{L^{2}},
\]
but with no control of the $L^{4}_{t}L^{\infty}_{x}$ norm one has to use local smoothing estimates close to the boundary, and Strichartz estimates for the usual Laplacian on $\mathbb{R}^{3}$ away from it. Here we give a simpler proof on the exterior of a strictly convex obstacle and for every $1+\frac{4}{3}<p<5$ using the Strichartz estimates \eqref{se}.

\subsection*{Acknowledgements} The author would like to thank Nicolas Burq for many helpful and stimulating discussions and Fabrice Planchon for having suggested to her the applications and for his careful reading of the paper. She would also like to thank Michael Taylor for having sent her the manuscript "Boundary problems  for the wave equations with grazing and gliding rays" and the referees for helpful comments and suggestions 
which greatly improved the presentation.

\section{Estimates for semi-classical Schr\"{o}dinger equation in a compact domain with strictly concave boundary}

In this section we prove Theorem \ref{thmstrichartz}. In what follows Assumptions \ref{assum1} are supposed to hold.
We may assume that the metric $g$ is extended smoothly across the boundary, so that $S$ is a geodesically concave subset of a complete, compact Riemannian manifold $\tilde{S}$. By the free semi-classical Schr\"odinger equation we mean the semi-classical Schr\"odinger equation on $\tilde{S}$, where the data $v_{0}$ has been extended to $\tilde{S}$ by an extension operator preserving the Sobolev spaces. By a broken geodesic in $S$ we mean a geodesic that is allowed to reflect off $\partial S$ according to the reflection law for the metric $g$.

\subsection{Restriction in a small neighborhood of the boundary.\\ Elliptic and hyperbolic regions}
We consider $\delta>0$ a small positive number and for $T>0$ small enough we set
\[
S(\delta,T):=\{(x,t)\in S\times [-T,T]| \text{dist}(x,\partial S)<\delta\}.
\]
On the complement of $S(\delta,T)$ in $S\times [-T,T]$, the solution $v(x,t)$ equals, in the semi-classical regime and modulo $O_{L^{2}}(h^{\infty})$ errors, the solution of the semi-classical Schr\"{o}dinger equation on a manifold without boundary for which sharp semi-classical Strichartz estimates follow by the work of Burq, G\'erard an Tzvetkow \cite{bgt04}, thus it suffices to establish Strichartz estimates for the norm of $v$ over $S(\delta,T)$.

We show that in order to prove Theorem \ref{thmstrichartz} it will be sufficient to consider only data $v_{0}$ supported outside a small neighborhood of the boundary. Recall that in \cite{le92} Lebeau proved that if $\Psi$ is supported in an interval $[\alpha_{0},\beta_{0}]$ and if $\varphi\in C^{\infty}_{0}(\mathbb{R})$ is equal to $1$ near the interval $[-\beta_{0},-\alpha_{0}]$, then for $t$ in a bounded set (and for $D_{t}=\frac{1}{i}\partial_{t}$) one has
\be\label{lochdt}
\forall N\geq 1,\quad \exists C_{N}>0\quad |(1-\varphi)(hD_{t})\exp{(ith\Delta_{g})}\Psi(-h^{2}\Delta_{g})v_{0}|\leq C_{N}h^{N}.
\ee
For $\delta$ and $T$ sufficiently small, let $\chi(x,t)\in C^{\infty}_{0}$ be compactly supported and be equal to $1$ on $S(\delta,T)$. Let $t_{0}>0$ be  such that $T=t_{0}/4$ and let $A\in C^{\infty}(\mathbb{R}^{n})$, $A=0$ near $\partial S$, $A=1$ outside a neighborhood of the boundary be such that every broken bicharacteristic $\gamma$ starting at $t=0$ from the support of $\chi(x,t)$ and for $-\tau\in [\alpha_{0},\beta_{0}]$, (where here $\tau$ denotes the dual time variable),  
satisfies
\be\label{condsupp}
\text{dist}(\gamma(t),\text{supp}(1-A))>0,\quad \forall t\in [-2t_{0},-t_{0}].
\ee
Let $\psi\in C^{\infty}(\mathbb{R})$, $\psi(t)=0$ for $t\leq-2t_{0}$, $\psi(t)=1$ for $t>-t_{0}$ and set 
\[
w(x,t)=\psi(t)\exp{(ith\Delta_{g})}\Psi(-h^{2}\Delta_{g})v_{0}.
\]
Then $w$ satisfies
\[
\left\{\begin{array}{ll}
ih\partial_{t}w+h^{2}\Delta_{g}w=ih\psi'(t)e^{ith\Delta_{g}}\Psi(-h^{2}\Delta_{g})v_{0},\\
w|_{\partial S\times\mathbb{R}}=0,\quad w|_{t\leq-2t_{0}}=0,
\end{array}
\right.
\]
and writing Duhamel formula we have
\[
w(x,t)=\int_{-2t_{0}}^{t}e^{i(t-s)h\Delta_{g}}\psi'(s)e^{ish\Delta_{g}}\Psi(-h^{2}\Delta_{g})v_{0}ds.
\]
Notice that $w(x,t)=v(x,t)$ if $t\geq -t_{0}$, hence for $t\in [-t_{0},T]$ we can write
\be\label{vtmic}
v(x,t)=\int_{-2t_{0}}^{-t_{0}}e^{i(t-s)h\Delta_{g}}\psi'(s)e^{ish\Delta_{g}}\Psi(-h^{2}\Delta_{g})v_{0}ds.
\ee
In particular, for $t\in [-T,T]$, $T=t_{0}/4$, $v(x,t)=w(x,t)$ is given by \eqref{vtmic}. We want to estimate the $L^{q}_{t}L^{r}_{x}$ norms of $v(x,t)$ for $(x,t)$ on $S(\delta,T)$ where $v=\chi v$. Let
\[
v_{Q}(x,t)=\int_{-2t_{0}}^{-t_{0}}e^{i(t-s)h\Delta_{g}}\psi'(s)Q(x)e^{ish\Delta_{g}}\Psi(-h^{2}\Delta_{g})v_{0}ds,\quad Q\in\{A,1-A\},
\]
then $v=v_{A}+ v_{1-A}$, where $v_{1-A}$ solves 
\[
\left\{\begin{array}{ll}
ih\partial_{t}v_{1-A}+h^{2}\Delta_{g}v_{1-A}=ih\psi'(t)(1-A)e^{ith\Delta_{g}}\Psi(-h^{2}\Delta_{g})v_{0},\\
v_{1-A}|_{\partial S\times\mathbb{R}}=0,\quad v_{1-A}|_{t<-2t_{0}}=0.
\end{array}
\right.
\]
We apply Proposition \ref{propschtdd} from the Appendix with $Q=1-A$, $\tilde{\psi}=\psi'$ to deduce that if $\rho_{0}\in WF_{b}(v_{1-A})$ then the broken bicharacteristic starting from $\rho_{0}$ must intersect the wave front set $WF_{b}((1-A)v)\cap \{t\in [-2t_{0},-t_{0}]\}$. Since we are interested in estimating the norm of $v$ on $S(\delta, T)$ it is enough to consider only $\rho_{0}\in WF_{b}(\chi v_{1-A})$. Thus, if $\gamma$ is a broken bicharacteristic starting at $t=0$ from $\rho_{0}$, $-\tau\in [\alpha_{0},\beta_{0}]$, then Proposition \ref{propschtdd} implies that for some $t\in [-2t_{0},-t_{0}]$, $\gamma(t)$ must intersect $WF_{b}((1-A)v)$. On the other hand from \eqref{condsupp} this implies (see Definition \ref{dfnwf}) that for every $\sigma\geq 0$
\be\label{est1a}
\forall N\geq 0 \quad \exists C_{N}>0\quad \|\chi v_{1-A}\|_{H^{\sigma}(S\times\mathbb{R})}\leq C_{N}h^{N}.
\ee
We are thus reduced to estimating $v(x,t)$ for initial data supported outside a small neighborhood of the boundary. Indeed, suppose that the estimates \eqref{sescl} hold true for any initial data compactly supported where $A\neq 0$. It follows from \eqref{vtmic}, \eqref{est1a} that
\begin{eqnarray*}
\|\chi v_{A}\|_{L^{q}((-T,T),L^{r}(S))} & \leq & 
\|\psi'(s)A(x)e^{ish\Delta_{g}}\Psi(-h^{2}\Delta_{g})v_{0}\|_{L^{1}(s\in(-2t_{0},-t_{0}),L^{2}(S))}\\ 
&  \lesssim &
(\int_{-2t_{0}}^{-t_{0}}|\Psi'(s)|ds)
\|\Psi(-h^{2}\Delta_{g})v_{0}\|_{L^{2}(S)}\\
 & = & \|\Psi(-h^{2}\Delta_{g})v_{0}\|_{L^{2}(S)},
\end{eqnarray*}
where we used the fact that the semi-classical Schr\"{o}dinger flow $\exp{(ihs\Delta_{g})}\Psi(-h^{2}\Delta_{g})$, which maps data at time $0$ to data at time $s$, is an isomorphisme on $H^{\sigma}(S)$ for every $\sigma\geq 0$.
\begin{rmq}
Notice that when dealing with the wave equation, since the speed of propagation is exact, one can take $\psi(t)=1_{\{t\geq -t_{0}\}}$ for some small $t_{0}\geq 0$ and reduce the problem to proving Strichartz estimates for the flow $\exp{(ih(t_{0}+.)\Delta_{g})}\Psi(-h^{2}\Delta_{g})$ and initial data compactly supported outside a small neighborhood of $\partial S$. This was precisely the strategy followed by Smith and Sogge in \cite{smso95}. 
\end{rmq}
Let $\Delta_{0}$ denote the Laplacian on $\tilde{S}$ coming from the extension of the metric $g$ smoothly across the boundary $\partial S$. 
We let $\mathcal{S}$ denote the outgoing solution to the Dirichlet problem for the semiclassical Schr\"{o}dinger operator on $S\times\mathbb{R}$. Thus, if $g$ is a function on $\partial S\times\mathbb{R}$ which vanishes for $t\leq -2t_{0}$, then $\mathcal{M}g$ is the solution on $S\times\mathbb{R}$ to
\be\label{scl2}
\left\{\begin{array}{c}
ih\partial_{t}\mathcal{M}g+h^{2}\Delta_{g}\mathcal{M}g=0,\\
\mathcal{M}g|_{\partial S\times\mathbb{R}}=g.
\end{array}
\right.
\ee
Then, for $t\in[-t_{0},T]$ and data $f$ supported outside a small neighborhood of the boundary and localized at frequency $1/h$, i.e. such that $f=\Psi(-h^{2}\Delta_{g})f$, we have 
\begin{multline*}
\chi v_{A}(x,t)=\chi\int_{-2t_{0}}^{-t_{0}}e^{i(t-s)h\Delta_{g}}\psi'(s)A(x)e^{ish\Delta_{g}}fds\\
=\chi\int_{-2t_{0}}^{-t_{0}}e^{i(t-s)h\Delta_{0}}\psi'(s)A(x)e^{ish\Delta_{0}}fds\\
 {}-\mathcal{M}\Big(\chi\int_{-2t_{0}}^{-t_{0}}e^{i(t-s)h\Delta_{0}}\psi'(s)A(x)e^{ish\Delta_{0}}fds|_{\partial S\times\mathbb{R}}\Big).
\end{multline*}
The cotangent bundle of $\partial S\times\mathbb{R}$ is divided into three disjoint sets: the hyperbolic and elliptic regions where the Dirichlet problem is respectively hyperbolic and elliptic, and the glancing region which is the boundary between the two. 

Let local coordinates be chosen such that $S=\{(x',x_{n})|x_{n}>0\}$ and $\Delta_{g}=\partial^{2}_{x_{n}}-r(x,D_{x'})$.
A point $(x',t,\eta',\tau)\in T^{*}(\partial S\times\mathbb{R})$ is classified as one of three distinct types. It is said to be \emph{hyperbolic} if $-\tau+r(x',0,\eta')>0$, so that there are two distinct nonzero real solutions $\eta_{n}$ to $\tau-r(x',0,\eta')=\eta^{2}_{n}$. These two solutions yield two distinct bicharacteristics, one of which enters $S$ as $t$ increases (the \emph{incoming ray}) and one which exits $S$ as $t$ increases (the \emph{outgoing ray}). The point is \emph{elliptic} if $-\tau+r(x',0,\eta')<0$, so there are no real solutions $\eta_{n}$ to $\tau-r(x',0,\eta')=\eta^{2}_{n}$. In the remaining case $-\tau+r(x',0,\eta')=0$, there is a unique solution which yields a glancing ray, and the point is said to be a \emph{glancing point}. 
We decompose the identity operator into
\[
\text{Id}(x,t)=\frac{1}{(2\pi h)^{n}}\int e^{\frac{i}{h}((x'-y')\eta'+(t-s)\tau)}(\chi_{h}+\chi_{e}+\chi_{gl})(y',\eta',\tau)d\eta' d\tau,
\]
where at $(y',\eta',\tau)$ we have
\[
\chi_{h}:=1_{\{-\tau+r(y',0,\eta')\geq c\}},\quad
\chi_{e}:=1_{\{-\tau+r(y',0,\eta')\leq -c\}},\quad
\chi_{h}:=1_{\{-\tau+r(y',0,\eta')\in [-c,c]\}},
\]
for some $c>0$ sufficiently small. The corresponding operators with symbols $\chi_{h}$, $\chi_{e}$, denoted $\Pi_{h}$, $\Pi_{e}$, respectively, are pseudo-differential cutoffs essentially supported inside the hyperbolic and elliptic regions, while the operator with symbol $\chi_{gl}$, denoted $\Pi_{gl}$, is essentially supported in a small set around the glancing region. 
Thus, on $S(\delta,T)$ we can write $\chi v_{A}$ as the sum of four terms 
\begin{multline}
\label{decompsolv}  
\chi\int_{-2t_{0}}^{-t_{0}}e^{i(t-s)h\Delta_{g}}\psi'(s)A(x)e^{ish\Delta_{g}}fds=\chi\int_{-2t_{0}}^{-t_{0}}e^{i(t-s)h\Delta_{0}}\psi'(s)A(x)e^{ish\Delta_{0}}fds\\
{}-\sum_{\Pi\in\{\Pi_{e},\Pi_{h},\Pi_{gl}\}}\mathcal{M}\Pi\Big(\chi\int_{-2t_{0}}^{-t_{0}}e^{i(t-s)h\Delta_{0}}\psi'(s)A(x)e^{ish\Delta_{0}}fds|_{\partial S\times\mathbb{R}}\Big).
\end{multline}
\begin{rmq}
For the first term in the right hand side, $\chi\int_{-2t_{0}}^{-t_{0}}e^{i(t-s)h\Delta_{0}}\psi'(s)A(x)e^{ish\Delta_{0}}fds$, the desired estimates follow as in the boundaryless case by the results of Staffilani and Tataru \cite{stta02} (since we considered the extension of the metric $g$ across the boundary to be smooth). 
\end{rmq}
\subsubsection{Elliptic region}
Using the compactness argument of the proof of Proposition \ref{prople92} from the Appendix, together with the inclusion (which follows from Proposition \ref{propellreg} in the Appendix)
\[
WF_{b}\Big(\chi\int_{-2t_{0}}^{-t_{0}}e^{i(t-s)h\Delta_{0}}\psi'(s)A(x)e^{ish\Delta_{0}}fds|_{\partial S\times\mathbb{R}}\Big)\subset \mathcal{H}\cup\mathcal{G},
\]
where $\mathcal{H}$ and $\mathcal{G}$ denote the hyperbolic and the glancing regions, respectively, it follows that the elliptic part satisfies for all $\sigma\geq 0$
\[
\mathcal{M}\Pi_{e}\Big(\chi\int_{-2t_{0}}^{-t_{0}}e^{i(t-s)h\Delta_{0}}\psi'(s)A(x)e^{ish\Delta_{0}}fds|_{\partial S\times\mathbb{R}}\Big)=O(h^{\infty})\|f\|_{H^{\sigma}(S)}.
\]
For the definition and properties of the $b$-wave front set see Appendix.

\subsubsection{Hyperbolic region}
If local coordinates are chosen such that $S=\{(x',x_{n})|x_{n}>0\}$,
on the essential support of $\Pi_{h}$ the forward Dirichlet problem can be solved locally, modulo smoothing kernels, on an open set in $\tilde{S}\times\mathbb{R}$ around $\partial S$. 
Precisely, microlocally near a hyperbolic point, the solution $v$ to \eqref{schrodscl} can be decomposed modulo smoothing operators into an incoming part $v_{-}$ and an outgoing part $v_{+}$ where
\[
v_{\pm}(x,t)=\frac{1}{(2\pi h)^{d}}\int e^{\frac{i}{h}\varphi_{\pm}(x,t,\xi)}\sigma_{\pm}(x,t,\xi,h)d\xi,
\]
where the phases $\varphi_{\pm}$ satisfy the eikonal equations 
\[
\left\{\begin{array}{ll}
\partial_{s}\varphi_{\pm}+<d\varphi_{\pm},d\varphi_{\pm}>_{g}=0,\\
\varphi_{+}|_{\partial S}=\varphi_{-}|_{\partial S},\quad \partial_{x_{n}}\varphi_{+}|_{\partial S}=-\partial_{x_{n}}\varphi_{-}|_{\partial S},
\end{array}
\right.
\]
where $<.,.>_{g}$ denotes the inner product induced by the metric $g$. The symbols are asymptotic expansions in $h$ and write $\sigma_{\pm}(.,h)=\sum_{k\geq 0}h^{k}\sigma_{\pm,k}$, where $\sigma_{0}$ solves the linear transport equation
\[
\partial_{s}\sigma_{\pm,0}+(\Delta_{g}\varphi_{\pm})\sigma_{\pm,0}+<d\varphi_{\pm},d\sigma_{\pm,0}>_{g}=0,
\]
while for $k\geq 1$, $\sigma_{\pm,k}$ satisfies the non-homogeneous transport equations
\[
\partial_{s}\sigma_{\pm,k}+(\Delta_{g}\varphi_{\pm})\sigma_{\pm,k}+<d\varphi_{\pm},d\sigma_{\pm,k}>_{g}=i\Delta_{g}\sigma_{\pm,k-1}.
\]
A direct computation shows that
\[
\|\sum_{\pm}v_{\pm}\|^{2}_{H^{\sigma}(S\times\mathbb{R})}\simeq \sum_{\pm}\|v_{\pm}\|^{2}_{H^{\sigma}(S\times\mathbb{R})}\simeq\|v\|^{2}_{H^{\sigma}(S\times\mathbb{R})}
\simeq\|v\|^{2}_{L^{\infty}(\mathbb{R})H^{\sigma}(S)}.
\]
Each component $v_{\pm}$ is a solution of linear Schr\"odinger equation (without boundary) and consequently satisfies the usual Strichartz estimates (see Burq, G\'erard and Tzvetkov \cite{bgt04}). 

Note that $\sum_{\pm}v_{\pm}$ contains the contribution from 
\[
\mathcal{M}\Pi_{h}\Big(\chi\int_{-2t_{0}}^{-t_{0}}e^{i(t-s)h\Delta_{0}}\psi'(s)A(x)e^{ish\Delta_{0}}\Psi(-h^{2}\Delta_{g})v_{0}ds|_{\partial S\times\mathbb{R}}\Big)
\] 
and a contribution from $\chi\int_{-2t_{0}}^{-t_{0}}e^{i(t-s)h\Delta_{0}}\psi'(s)A(x)e^{ish\Delta_{0}}\Psi(-h^{2}\Delta_{g})v_{0}ds$.

\subsection{Glancing region}

Near a diffractive point we use the Melrose and Taylor construction for the wave equation in order to write, following Zworski \cite{zw90}, the solution to the wave equation as a finite sum of pseudo-differential cutoffs, each essentially supported in a suitably small neighborhood of a glancing ray. Using the Fourier transform in time we obtain a parametrix for the semi-classical Schr\"{o}dinger equation \eqref{schrodscl} microllocally near a glancing direction and modulo smoothing operators.

\subsubsection{Preliminaries. Parametrix for the wave equation near the glancing region}
We start by recalling the results by Melrose and Taylor \cite{meltay85}, \cite{meta85}, Zworski \cite[Prop.4.1]{zw90} for the wave equation near the glancing region.
Let $w$ solve the (semi-classical) wave equation on $S$ with Dirichlet boundary conditions
\begin{equation}\label{unde}
\left\{\begin{array}{ll}
h^{2}D_{t}^{2}w+h^{2}\Delta_{g}w=0,\quad S\times\mathbb{R},\quad w|_{\partial S\times\mathbb{R}}=0,\\
w(x,0)=f(x),\quad D_{t}w(x,0)=g(x),
\end{array}
\right.
\end{equation}
where $f$, $g$ are compactly supported in $S$ and localized at spacial frequency $1/h$, and where $D_{t}=\frac{1}{i}\partial_{t}$.
\begin{prop}\label{propglan}
Microlocally near a glancing direction the solution to \eqref{unde} writes, modulo smoothing operators
\begin{multline}\label{parametrix}
w(x,t)=\frac{1}{(2\pi h)^{n}}\int_{\mathbb{R}^{n}}e^{\frac{i}{h}(\theta(x,\xi)+it\xi_{1})}\Big[a(x,\xi/h)\Big(A_{-}(\zeta(x,\xi/h))-A_{+}(\zeta(x,\xi/h))\frac{A_{-}(\zeta_{0}(\xi/h))}{A_{+}(\zeta_{0}(\xi/h))}\Big)\\
+b(x,\xi/h)\Big(A'_{-}(\zeta(x,\xi/h))-A'_{+}(\zeta(x,\xi/h))\frac{A_{-}(\zeta_{0}(\xi/h))}{A_{+}(\zeta_{0}(\xi/h))}\Big)\Big]\times
\widehat{K(f,g)}(\frac{\xi}{h})d\xi,
\end{multline}
where the symbols $a$, $b$ and the phases $\theta$, $\zeta$ have the following properties: $a$ and $b$ are symbols of type $(1,0)$ and order $1/6$ and $-1/6$, respectively, both of which are supported in a small conic neighborhood of the $\xi_{1}$ axis and where $K$ is a classical Fourier integral operator of order $0$ in $f$ and order $-1$ in $g$, compactly supported on both sides. The phases $\theta$ and $\zeta$ are real, smooth and homogeneous of degree $1$ and $2/3$, respectively.  If we denote $Ai(z)$ the Airy function, then $A_{\pm}$ are defined by $A_{\pm}(z)=Ai(e^{\mp 2\pi i/3}z)$.
\end{prop}
\begin{rmq}
If local coordinates are chosen so that $\Omega$ is given by $x_{n}>0$, the phase functions $\theta$, $\zeta$ satisfy the eikonal equations
\be\label{eikonal}
\left\{ \begin{array}{c}
\xi^{2}_{1}-<d\theta,d\theta>_{g}+\zeta<d\zeta,d\zeta>_{g}=0,\\
<d\theta,d\zeta>_{g}=0,\\
\zeta(x',0,\xi)=\zeta_{0}(\xi)=-\xi^{-1/3}_{1}\xi_{n},
\end{array}
\right. 
\ee
in the region $\zeta\leq 0$. Here $x'=(x_{1},..,x_{n-1})$ and $<.,.>_{g}$ denotes the inner product given by the metric $g$. The phases also satisfy the eikonal equations \eqref{eikonal} to infinite order at $x_{n}=0$ in the region $\zeta>0$. 
\end{rmq}
\begin{rmq}
Notice that one can think of $A_{-}(\zeta)$ (at least away from the boundary $x_{n}=0$) as the incoming contribution and of $A_{+}(\zeta)\frac{A_{-}(\zeta_{0})}{A_{+}(\zeta_{0})}$ as the outgoing one. From \cite[Section 2]{zw90} we have
\[
\frac{A_{-}}{A_{+}}(z)\simeq 
\left\{ \begin{array}{c}
-e^{i\pi /3}+O(z^{-\infty}), \quad z\rightarrow \infty\\
e^{i(4/3)(-z)^{3/2}}\sum_{j\geq 0}\beta_j z^{-3j/2}, \quad z\rightarrow -\infty,
\end{array}
\right.
\]
where the part $z\rightarrow\infty$ corresponds to the free wave, while the oscillatory one to the billiard ball map shift corresponding to reflection.
Using $Ai(\zeta)=e^{i\pi/3}A_{+}(\zeta)+e^{-i\pi/3}A_{-}(\zeta)$, we write
\[
A_{-}(\zeta)-A_{+}(\zeta)\frac{A_{-}(\zeta_{0})}{A_{+}(\zeta_{0})}=e^{i\pi/3}\Big(Ai(\zeta)-A_{+}(\zeta)\frac{Ai(\zeta_{0})}{A_{+}(\zeta_{0})}\Big).
\]
\end{rmq}

\subsubsection{Parametrix for the solution to the semi-classical Schr\"odinger equation near the glancing region}
Let now $v(x,t)$ be the solution of the semi-classical Schr\"odinger equation \eqref{schrodscl} where the initial data $v_{0}\in L^{2}(S)$ is spectrally localized at spatial frequency $1/h$, i.e. $v_{0}(x)=\Psi(-h^{2}\Delta_{g})v_{0}(x)$. From the discussion at the beginning of this section we see that it will be enough to consider $v_{0}$ compactly supported outside some small neighborhood of $\partial S$. Under this assumption $\Psi(-h^{2}\Delta_{g})v_{0}$ is a well-defined pseudo-differential operator for which the results of Burq, G\'erard and Tzvetkov \cite[Section 2]{bgt04} apply.

Let $(e_{\lambda}(x))_{\lambda\geq 0}$ be the eigenbasis of $L^{2}(S)$ consisting in eigenfunctions of $-\Delta_{g}$ associated to the eigenvalues $(\lambda^{2})$, so that $-\Delta_{g}e_{\lambda}=\lambda^{2}e_{\lambda}$. We write
\be
\Psi(-h^{2}\Delta_{g})v_{0}(x)= \sum_{h^{2}\lambda^{2}\in[\alpha_{0},\beta_{0}]}\Psi(h^{2}\lambda^{2})v_{\lambda}e_{\lambda}(x),
\ee
and hence
\be
e^{ith\Delta_{g}}\Psi(-h^{2}\Delta_{g})v_{0}(x)= \sum_{h^{2}\lambda^{2}\in[\alpha_{0},\beta_{0}]}\Psi(h^{2}\lambda^{2})e^{-ith\lambda^{2}}v_{\lambda}e_{\lambda}(x).
\ee
If $\delta$ denotes the Dirac function, then the Fourier transform of $v(x,t)$ writes
\be
\hat{v}(x,\frac{\tau}{h})= h\sum_{h^{2}\lambda^{2}\in[\alpha_{0},\beta_{0}]}\Psi(h^{2}\lambda^{2})\delta_{\{-\tau=h^{2}\lambda^{2}\}}v_{\lambda}e_{\lambda}(x).
\ee
For $t\in\mathbb{R}$ we can define (since $\hat{v}$ has compact support away from $0$)
\begin{align}\label{defwondes}
\nonumber
w(x,t): & =\frac{1}{2\pi h}\int_{0}^{\infty}e^{\frac{it\sigma}{h}}\hat{v}(x,-\frac{\sigma^{2}}{h})d\sigma\\
\nonumber
  & =-\frac{1}{4\pi h}\int_{-\infty}^{0} e^{\frac{it\sqrt{-\tau}}{h}}\frac{1}{\sqrt{-\tau}}\hat{v}(x,\frac{\tau}{h})d\tau\\
\nonumber
& =  -\frac{1}{2}\sum_{h^{2}\lambda^{2}\in[\alpha_{0},\beta_{0}]}\Psi(h^{2}\lambda^{2})\Big(\frac{1}{2\pi}\int_{-\infty}^{0}e^{i\frac{t\sqrt{-\tau}}{h}}\frac{1}{\sqrt{-\tau}}\delta_{\{-\tau=h^{2}\lambda^{2}\}}d\tau \Big)v_{\lambda}e_{\lambda}(x)\\
&= -\frac{1}{2}\sum_{h^{2}\lambda^{2}\in[\alpha_{0},\beta_{0}]}\frac{1}{h\lambda}\Psi(h^{2}\lambda^{2})e^{it\lambda}v_{\lambda}e_{\lambda}(x).
\end{align}
Then $w(x,t)$ solves the wave equation
\be\label{undedoi}
\left\{\begin{array}{ll}
h^{2}D_{t}^{2}w+h^{2}\Delta_{g}w=0,\quad\text{on}\quad S\times\mathbb{R},\quad w|_{\partial S\times\mathbb{R}}=0,\\
w(x,0)=f_{h}(x),\quad D_{t}w(x,0)=g_{h}(x),
\end{array}
\right.
\ee
where the initial data are given by
\be\label{deffh}
f_{h}(x)= -\frac{1}{2}\sum_{h^{2}\lambda^{2}\in[\alpha_{0},\beta_{0}]}\frac{1}{h\lambda}\Psi(h^{2}\lambda^{2})v_{\lambda}e_{\lambda}(x),
\ee
\be\label{defgh}
g_{h}(x)=-\frac{1}{2h}\sum_{h^{2}\lambda^{2}\in[\alpha_{0},\beta_{0}]}\Psi(h^{2}\lambda^{2})v_{\lambda}e_{\lambda}(x)= -\frac{1}{2h}\Psi(-h^{2}\Delta_{g})v_{0}(x).
\ee
From \eqref{deffh}, \eqref{defgh} it follows that
\be\label{equivnormes}
h\|g_{h}\|_{L^{2}(S)}\simeq\|f_{h}\|_{L^{2}(S)}\simeq \|\Psi(-h^{2}\Delta_{g})v_{0}\|_{L^{2}(S)},
\ee
where by $\alpha\simeq\beta$ we mean that there is $C>0$ such that $C^{-1}\alpha<\beta<C\alpha$.

Indeed, in order to prove \eqref{equivnormes} notice that $w$ defined by \eqref{defwondes} satisfies in fact
\[
(hD_{t}-h\sqrt{-\Delta_{g}})w=0
\]
and (since $\Delta_{g}$ and $D_{t}$ commute) we have
\[
f_{h}=w|_{t=0}=[(\sqrt{-\Delta_{g}})^{-1}D_{t}w]|_{t=0}=(\sqrt{-\Delta_{g}})^{-1}(D_{t}w|_{t=0})=(\sqrt{-\Delta_{g}})^{-1}g_{h}.
\]
Due to the spectral localization and since $g_{h}=-\frac{1}{2h}\Psi(-h^{2}\Delta_{g})v_{0}$ we deduce \eqref{equivnormes}.

By the $L^{2}$ continuity of the (classical) Fourier integral operator $K$ introduced in Proposition \ref{propglan} we deduce
\be\label{kl2}
\|K(f_{h},g_{h})\|_{L^{2}(S)}\leq C(\|f_{h}\|_{L^{2}(S)}+h\|g_{h}\|_{L^{2}(S)})\simeq \|\Psi(-h^{2}\Delta_{g})v_{0}\|_{L^{2}(S)}.
\ee
The solution $v(x,t)$ of \eqref{schrodscl} writes
\begin{align}
\nonumber
v(x,t) & =\frac{1}{2\pi h}\int_{0}^{\infty} e^{-\frac{it\sigma^{2}}{h}}2\sigma\hat{v}(x,-\frac{\sigma^{2}}{h})d\sigma\\
& = \frac{1}{2\pi h}\int_{0}^{\infty} e^{-i\frac{t\sigma^{2}}{h}}2\sigma \int_{s\in\mathbb{R}}e^{-i\frac{s\sigma}{h}}w(x,s)dsd\sigma.
\end{align}
The next step is to use Proposition \ref{unde} in order to obtain a representation of $v(x,t)$ near the glancing region: notice that the glancing part of the stationary wave $\hat{w}(x,\frac{\sigma}{h})$ is given by
\begin{align}\label{eqwvhat}
\nonumber
1_{\{\sigma^{2}+r(x',0,\eta')\in [-c,c]\}}\hat{w}(x,\frac{\sigma}{h}) & =1_{\{\sigma^{2}+r(x',0,\eta')\in [-c,c]\}}\hat{v}(x,-\frac{\sigma^{2}}{h})\\
& =1_{\{-\tau+r(x',0,\eta')\in [-c,c]\}}\hat{v}(x,\frac{\tau}{h}),
\end{align}
with $\tau=-\sigma^2$ and where $c>0$ is sufficiently small. The equality in \eqref{eqwvhat} follows from \eqref{defwondes} and from the fact that $\hat{v}$ is essentially supported for the second variable in the interval $[-\beta_{0},-\alpha_{0}]$. Consequently we can apply Proposition \ref{unde} and determine a representation for $v$ near the glancing region (for the Schr\"odinger equation) as follows
\begin{multline}\label{vglancing}
v(x,t)=
\frac{1}{(2\pi h)^{n}}\int_{\mathbb{R}^{n}}e^{\frac{i}{h}(\theta(x,\xi)-t\xi^{2}_{1})}2\xi_{1}\Big[a(x,\xi/h)\Big(Ai(\zeta(x,\xi/h))-A_{+}(\zeta(x,\xi/h))\frac{Ai(\zeta_{0}(\xi/h))}{A_{+}(\zeta_{0}(\xi/h))}\Big)\\+
b(x,\xi/h)\Big( Ai'(\zeta(x,\xi/h))-A'_{+}(\zeta(x,\xi/h))\frac{Ai(\zeta_{0}(\xi/h))}{A_{+}(\zeta_{0}(\xi/h))}\Big)\Big]\widehat{K(f_{h},g_{h})}(\frac{\xi}{h})d\xi,
\end{multline}
where $a$, $b$ and $K$ are those defined in Proposition \ref{propglan} and $f_{h}$, $g_{h}$ are given by \eqref{deffh}, \eqref{defgh}. The initial data $f_{h}$, $g_{h}$ are both supported, like $v_{0}$, away from $\partial S$, and consequently their $\dot{H}^{\sigma}(S)$ norms for $\alpha<n/2$ will be comparable to the norms of the non-homogeneous Sobolev space $H^{\sigma}(\mathbb{R}^{n})$, so we shall henceforth work with the latter norms on the data $f_{h}$, $g_{h}$.
\begin{rmq}
Notice that it is enough to prove semi-classical Strichartz estimates only for the "outgoing" piece  corresponding to the oscillatory term $A_{+}(\zeta)\frac{Ai(\zeta_{0})}{A_{+}(\zeta_{0})}$, since the direct term (corresponding to $Ai(\zeta)$) has already been dealt with (see the remark following \eqref{decompsolv}).
\end{rmq}
We deduce from \eqref{kl2}, \eqref{vglancing} that in order to finish the proof of Theorem \ref{thmstrichartz} we need only to show that the operator $A_{h}$ defined, for $f$ supported away from $\partial S$ and spectrally localized at frequency $1/h$, i.e. such that $f=\Psi(-h^{2}\Delta_{g})f$, by
\begin{multline}\label{opa}
A_{h}f(x,t)=\frac{1}{(2\pi h)^{n}}\int_{\mathbb{R}^{n}}
2\xi_{1}(a(x,\xi/h)A_{+}(\zeta(x,\xi/h))+b(x,\xi/h)A'_{+}(\zeta(x,\xi/h)))\\
\times e^{\frac{i}{h}(\theta(x,\xi)-t\xi^{2}_{1})}
\frac{Ai(\zeta_{0}(\xi/h))}{A_{+}(\zeta_{0}(\xi/h))}\widehat{f}(\frac{\xi}{h})d\xi,
\end{multline}
satisfies
\begin{equation}\label{stra}
\|A_{h}f\|_{L^{q}((0,T], L^{r}(\mathbb{R}^{n}))}\leq C h^{-\frac{1}{q}}\|f\|_{L^{2}(\mathbb{R}^{n})}.
\end{equation}
\begin{rmq}\label{rmqlocaliz}
We introduce a cut-off function $\chi_{1}\in C^{\infty}_{0}(\mathbb{R}^{n})$ equal to $1$ on the support of $f$ and equal to $0$ near $\partial S$. Since $\chi_{1}$ is supported away from the boundary it follows from \cite[Prop.2.1]{bgt04} (which applies here in its adjoint form) that $\Psi(-h^{2}\Delta_{g})\chi_{1}f$ is a pseudo-differential operator and writes (in a patch of local coordinates)
\be
 \Psi(-h^{2}\Delta_{g})\chi_{1}f=d(x,hD_{x})\chi_{2}f+O_{L^{2}(S)}(h^{\infty}),
\ee
where $\chi_{2}\in C^{\infty}_{0}(\mathbb{R}^{n})$ is equal to $1$ on the support of $\chi_{1}$ and where $d(x,D_{x})$ is defined for $x$ in the suitable coordinate patch using the usual pseudo-differential quantization rule, 
\[
d(x,D_{x})f(x)=\frac{1}{(2\pi)^{n}}\int_{\mathbb{R}^{n}}e^{ix\xi}d(x,\xi)\hat{f}(\xi)d\xi, \quad d\in C^{\infty}_{0},
\]
with symbol $d$ compactly supported for $|\xi|^{2}_{g}:=<\xi,\xi>_{g}\in [\alpha_{0},\beta_{0}]$, which follows by the condition of the support of $\Psi$. Since the principal part of the Laplace operator $\Delta_{g}$ is uniformly elliptic, we can introduce a smooth radial function $\psi\in C^{\infty}_{0}([\frac{1}{\delta}\alpha^{1/2}_{0},\delta\beta^{1/2}_{0}])$ for some $\delta\geq1$ such that
$\psi(|\xi|)d=d$ everywhere. In what follows we shall prove \eqref{stra} where, instead of $f$ we shall write $\psi(|\xi|)f$, keeping in mind that $f$ is supported away from the boundary and localized at spatial frequency $1/h$.
\end{rmq}
The proof of Theorem \ref{thmstrichartz} will be completed once we prove \eqref{stra}.
In order to do that, we split the operator $A_{h}$ into two parts: a main term and a diffractive term. To this end, let $\chi(s)$ be a smooth function satisfying 
\[
\text{supp}\chi\subset (-\infty, -1], \quad \text{supp} (1-\chi)\subset[-2,\infty).
\]
We write this operator as a sum $A_{h}=M_{h}+D_{h}$, by decomposing
\[
A_{+}(\zeta(x,\xi))=(\chi A_{+})(\zeta(x,\xi))+((1-\chi) A_{+})(\zeta(x,\xi)),
\]
and letting the "main term" be defined for $f$ like in Remark \ref{rmqlocaliz} by 
\begin{multline}
M_{h}f(x,t)=\frac{1}{(2\pi h)^{n}}\int_{\mathbb{R}^{n}}
2\xi_{1}(a(x,\xi/h)(\chi A_{+})(\zeta(x,\xi/h))+b(x,\xi/h)(\chi A'_{+})(\zeta(x,\xi/h)))\\
\times e^{\frac{i}{h}(\theta(x,\xi)-t\xi^{2}_{1})}
\frac{Ai(\zeta_{0}(\xi/h))}{A_{+}(\zeta_{0}(\xi/h))}\psi(|\xi|)\hat{f}(\frac{\xi}{h})d\xi.
\end{multline}
The "diffractive term" is then defined for $f$ like before by
\begin{multline}
D_{h}f(x,t)=\frac{1}{(2\pi h)^{n}}\int_{\mathbb{R}^{n}}
2\xi_{1}(a(x,\xi/h)((1-\chi) A_{+})(\zeta(x,\xi/h))+b(x,\xi/h)((1-\chi) A'_{+})(\zeta(x,\xi/h)))\\
\times e^{\frac{i}{h}(\theta(x,\xi)-t\xi^{2}_{1})}
\frac{Ai(\zeta_{0}(\xi/h))}{A_{+}(\zeta_{0}(\xi/h))}\psi(|\xi|)\hat{f}(\frac{\xi}{h})d\xi.
\end{multline}
We analyze these operators separately, following the ideas of \cite{smso95}:

\subsubsection{The main term $M_{h}$}\label{sectmainterm}
To estimate the "main term" $M_{h}$ we first use the fact that
\be\label{usor}
|\frac{Ai(s)}{A_{+}(s)}|\leq 2, \quad s\in\mathbb{R}.
\ee
Consequently, since the term $\frac{Ai(\zeta_{0})}{A_{+}(\zeta_{0})}$ acts like a multiplier and so does $\xi_{1}$ which is localized in the interval $[\alpha_{0},\beta_{0}]$ (this follows from \eqref{lochdt}), the estimates for $M_{h}$ will follow from showing that the operator 
\begin{multline}\label{formw}
f\rightarrow \frac{1}{(2\pi h)^{n}}\int_{\mathbb{R}^{n}}
(a(x,\xi/h)(\chi A_{+})(\zeta(x,\xi/h))+b(x,\xi/h)(\chi A'_{+})(\zeta(x,\xi/h)))\\
\times e^{\frac{i}{h}(\theta(x,\xi)-t\xi^{2}_{1})}\psi(|\xi|)\hat{f}(\frac{\xi}{h})d\xi
\end{multline}
satisfies the same bounds like in \eqref{stra} for $f$ spectrally localized at frequency $1/h$. Following \cite[Lemma 4.1]{zw90}, we write $\chi A_{+}$ and $(\chi A_{+})'$ in terms of their Fourier transform to express the phase function of this operator
\begin{equation}\label{phi}
\phi(t,x,\xi)=-t\xi^{2}_{1}+\theta(x,\xi)-\frac{2}{3}(-\zeta)^{3/2}(x,\xi),
\end{equation}
which satisfies the eikonal equation \eqref{eikonal}. We denote its symbol $c_{m}(x,\xi/h)$, $c_{m}(x,\xi)\in \mathcal{S}^{0}_{2/3,1/3}(\mathbb{R}^{n}\times\mathbb{R}^{n})$ and we also denote the operator defined in \eqref{formw} by $W^{m}_{h}$, thus 
\[
W^{m}_{h}f(x,t)=\frac{1}{(2\pi h)^{n}}\int_{\mathbb{R}^{n}}e^{\frac{i}{h}\phi(t,x,\xi)}c_{m}(x,\xi/h)\psi(|\xi|)\hat{f}(\frac{\xi}{h})d\xi.
\]
\begin{prop}\label{propestim}
Let $(q,r)$ be an admissible pair with $q>2$, let $T>0$ be sufficiently small and for $f=d(x,D_{x})\chi_{2}f+O_{L^{2}(\Omega)}(h^{\infty})$ like in Remark \ref{rmqlocaliz} let 
\[
W_{h}f(x,t):=W_{h}^{m}f(x,t)=\frac{1}{(2\pi h)^{n}}\int e^{\frac{i}{h}\phi(t,x,\xi)}c_{m}(x,\xi/h)\psi(|\xi|)\hat{f}(\frac{\xi}{h})d\xi.
\]
Then the following estimates hold
\be\label{eststri}
\|W_{h}f\|_{L^{q}((0,T], L^{r}(\mathbb{R}^{n}))}\leq C h^{-\frac{1}{q}}\|f\|_{L^{2}(\mathbb{R}^{n})}.
\ee
\end{prop}
In the rest of this section we prove Proposition \ref{propestim}. The first step in the proof is a $\text{TT}^{*}$ argument. Explicitly,
\[
\widehat{W^{*}_{h}(F)}(\frac{\xi}{h})=\int e^{-\frac{i}{h}\phi(s,y,\xi)}F(y,s)\overline{c_{m}(y,\xi/h)}dyds,
\]
and if we set
\begin{multline}
(T_{h}F)(x,t)=(W_{h}W^{*}_{h}F)(x,t)\\ =
\frac{1}{(2\pi h)^{n}}\int e^{\frac{i}{h}(\phi(t,x,\xi)-\phi(s,y,\xi))}c_{m}(x,\xi/h)\overline{c_{m}(y,\xi/h)}\psi^{2}(|\xi|)F(y,s)d\xi dsdy,
\end{multline}
then inequality \eqref{eststri} is equivalent to
\be\label{estttstar}
\|T_{h}F\|_{L^{q}((0,T], L^{r}(\mathbb{R}^{n}))}\leq C h^{-\frac{2}{q}}\|F\|_{L^{q'}((0,T], L^{r'}(\mathbb{R}^{n}))},
\ee
where $(q',r')$ satisfies $1/q+1/q'=1$, $1/r+1/r'=1$.
To see, for instance, that \eqref{estttstar} implies \eqref{eststri}, notice that the dual version of \eqref{eststri} is
\[
\|W^{*}_{h}F\|_{L^{2}(\mathbb{R}^{n})} \leq C h^{-\frac{1}{q}}\|F\|_{L^{q'}((0,T], L^{r'}(\mathbb{R}^{n}))},
\]
and we have
\begin{align}
\nonumber
\|W^{*}_{h}F\|^{2}_{L^{2}(\mathbb{R}^{n})} & =\int W_{h}W^{*}_{h}F\bar{F}dtdx\\
& \leq \|T_{h}F\|_{L^{q}((0,T], L^{r}(\mathbb{R}^{n}))}\|F\|_{L^{q'}((0,T], L^{r'}(\mathbb{R}^{n}))}.
\end{align}
Therefore we only need to prove \eqref{estttstar}.
Since the symbols are of type $(2/3,1/3)$ and not of type $(1,0)$, before starting the proof of \eqref{estttstar} for the operator $T_{h}$ we need to make a further decomposition: let $\rho\in C^{\infty}_{0}(\mathbb{R})$ satisfying $\rho(s)=1$ near $0$, $\rho(s)=0$ if $|s|\geq 1$. Let 
\[
T_{h}F=T^{f}_{h}F+T_{h}^{s}F,
\]
where 
\be
T_{h}^{s}F(x,t)=\int K_{h}^{s}(t,x,s,y)F(y,s)ds dy,
\ee
\begin{multline}
K_{h}^{s}(t,x,s,y)=\frac{1}{(2\pi h)^{n}}\int e^{\frac{i}{h}(\phi(t,x,\xi)-\phi(s,y,\xi))}(1-\rho(h^{-1/3}|t-s|))\\ \times c_{m}(x,\xi/h)\overline{c_{m}(y,\xi/h)}\psi^{2}(|\xi|)d\xi,
\end{multline}
while 
\be
T_{h}^{f}F(x,t)=\int K_{h}^{f}(t,x,s,y)F(y,s)ds dy,
\ee
\begin{multline}
K_{h}^{f}(t,x,s,y)=\frac{1}{(2\pi h)^{n}}\int e^{\frac{i}{h}(\phi(t,x,\xi)-\phi(s,y,\xi))}\rho(h^{-1/3}|t-s|)\\ \times c_{m}(x,\xi/h)\overline{c_{m}(y,\xi/h)}\psi^{2}(|\xi|)d\xi.
\end{multline}
\begin{rmq}
The two pieces will be handled differently. The kernel of $T^{f}_{h}$ is supported in a suitable small set and it will be estimate by "freezing"  the coefficients. To estimate $T^{s}_{h}$ we shall use the stationary phase method for type $(1,0)$ symbols. For type $(2/3,1/3)$ symbols, these stationary phase arguments break down if $|t-s|$ is smaller than $h^{1/3}$, which motivates the decomposition. We use here the same arguments as in \cite{smso95}. 
\end{rmq}

\begin{itemize}
\item The "stationary phase admissible" term $T^{s}_{h}$
\begin{prop}\label{propfio}
There is a constant $1<C_{0}<\infty$ such that the kernel $K^{s}_{h}$ of $T^{s}_{h}$ satisfies
\be\label{ippo}
|K^{s}_{h}(t,x,s,y)|\leq C_{N}h^{N}\quad \forall N, \quad \text{if}\quad  \frac{|t-s|}{|x-y|} \notin [C_{0}^{-1},C_{0}].
\ee
Moreover, there is a function $\xi_{c}(t,x,s,y)$ which is smooth in the variables $(t,s)$, uniformly over $(x,y)$, so that if $C_{0}^{-1}\leq \frac{|t-s|}{|x-y|}\leq C_{0}$, then 
\be
|K^{s}_{h}(t,x,s,y)| \lesssim h^{-n}(1+\frac{|t-s|}{h})^{-n/2}, \quad \text{for}\quad |t-s|\geq h^{1/3}.
\ee
\end{prop}
\begin{proof} 
We shall use stationary phase lemma to evaluate the kernel $K^{s}_{h}$ of $T^{s}_{h}$. The critical points occur when $|t-s|\simeq |x-y|$. For some constant $C_{0}$ and for $|\xi|\in \text{supp}\psi$, $\xi_{1}$ in a small neighborhood of $1$, we have
\[
|\nabla_{\xi}(\phi(t,x,\xi)-\phi(s,y,\xi))|\simeq |t-s|+|x-y|\geq h^{1/3}, \quad\text{if}\quad \frac{|t-s|}{|x-y|} \notin [C_{0}^{-1},C_{0}].
\]  
Since $c\in S^{0}_{2/3,1/3}$, an integration by parts leads to \eqref{ippo}. If $|t-s|\simeq |x-y|$ we introduce a cutoff function $\kappa(\frac{|x-y|}{|t-s|})$, $\kappa\in C^{\infty}_{0}(\mathbb{R}\setminus\{0\})$. The phase function can be written as 
\[
\phi(t,x,\xi)-\phi(s,y,\xi)=(t-s)\Theta(t,x,s,y,\xi)\quad\text{for} \quad |t-s|\simeq |x-y|\geq h^{1/3}.
\]
We want to apply the stationary phase method with parameter $|t-s|/h\geq h^{-2/3}\gg 1$ to estimate $K^{s}_{h}$. For $x$, $y$, $t$, $s$ fixed we must show that the critical points of $\Theta$ are non-degenerate. 
\begin{lemma}\label{lemptcrit}
If $T$ is sufficiently small then the phase function $\Theta(t,x,s,y,\xi)$ admits a unique, non-degenerate critical point $\xi_{c}$. Moreover,  for $0\leq t,s\leq T$, the function $\xi_{c}(t,x,s,y)$ solving $\nabla_{\xi}\Theta(t,x,s,y,\xi_{c})=0$ is smooth in $t$ and $s$, with uniform bounds on derivatives as $x$ and $y$ vary and we have
\begin{equation}\label{estimxic}
|\partial^{\alpha}_{t,s}\partial^{\gamma}_{x,y}\xi_{c}(t,x,s,y)|\leq C_{\alpha,\gamma} h^{-|\alpha|/3} \quad \text{if} \quad |x-y|\geq h^{1/3}.
\end{equation}
\end{lemma}
\begin{proof}
The phase $\Theta(t,x,s,y,\xi)$ writes
\begin{align}
\nonumber
\Theta(t,x,s,y,\xi)& =\xi^{2}_{1}+\frac{1}{(t-s)}(\phi(0,x,\xi)-\phi(0,y,\xi))\\
& =\xi^{2}_{1}+\frac{1}{(t-s)}\sum_{j=1}^{n}(x_{j}-y_{j})\partial_{x_{j}}\phi(0,z_{x,y},\xi),
\end{align}
for some $z_{x,y}$ close to $x$, $y$ (if $T$ is sufficiently small then $|t-s|\simeq |x-y|$ is small), and using the eikonal equations \eqref{eikonal} we can write
\[
\Theta(t,x,s,y,\xi)=<\nabla_{x}\phi,\nabla_{x}\phi>_{g}(0,z_{x,y},\xi)-\frac{1}{(t-s)}\sum_{j=1}^{n}(x_{j}-y_{j})\partial_{x_{j}}\phi(0,z_{x,y},\xi).
\]
Let us write $<\nabla_{x}\phi,\nabla_{x}\phi>_{g}=\sum_{j,k}g^{j,k}\partial_{x_{j}}\phi\partial_{x_{k}}\phi$ and compute explicitly $\nabla_{\xi}\Theta$. For each $l\in\{1,..,n\}$ we have
\begin{equation}\label{crit}
\partial_{\xi_{l}}\Theta(t,x,s,y,\xi)=\sum_{j=1}^{n}\partial^{2}_{\xi_{l},x_{j}}\phi(0,z_{x,y},\xi)\Big(2\sum_{k=1}^{n}g^{j,k}(z_{x,y})\partial_{x_{k}}\phi(0,z_{x,y},\xi)-\frac{(x_{j}-y_{j})}{(t-s)}\Big),
\end{equation}
thus
\begin{equation}\label{nablaThet}
\nabla_{\xi}\Theta(t,x,s,y,\xi)=\nabla^{2}_{\xi,x}\phi(0,z_{x,y},\xi).
\left(
\begin{array}{c}
2 \sum_{k}g^{1,k}(z_{x,y})\partial_{x_{k}}\phi(0,z_{x,y},\xi)-\frac{x_{1}-y_{1}}{(t-s)}    \\
  .  \\
  .  \\
 2 \sum_{k}g^{n,k}(z_{x,y})\partial_{x_{k}}\phi(0,z_{x,y},\xi)-\frac{x_{n}-y_{n}}{(t-s)} 
\end{array}
\right),
\end{equation}
where $\nabla^{2}_{\xi,x}\phi=(\partial^{2}_{\xi_{l},x_{j}}\phi)_{l,j\in\{1,..,n\}}$ is the matrix $n\times n$ whose elements are the second derivatives of $\phi$ with respect to $\xi$ and $x$. We need the next lemma:
\begin{lemma}\label{lemdetermin}(see \cite[Lemma 3.9]{smso94})
For $\xi$ in a conic neighborhood of the $\xi_{1}$ axis the mapping
\[
x\rightarrow \nabla_{\xi}\Big(\theta(x,\xi)-\frac{2}{3}(-\zeta)^{3/2}(x,\xi)\Big)
\]
is a diffeomorphisme on the complement of the hypersurface $\zeta=0$, with uniform bounds of the Jacobian of the inverse mapping. 
\end{lemma}
A direct corollary of Lemma \ref{lemdetermin} is the following:
\begin{cor}\label{cordetermin}
If $T$ is small enough and $|x-y|\simeq |t-s|\leq 2T$ then 
\be\label{determinhessphi}
\det(\nabla^{2}_{\xi,x}\phi)(0,z_{x,y},\xi)\neq 0.
\ee
\end{cor}
In what follows we complete the proof of Lemma \ref{lemptcrit}.
A critical point for $\Theta$ satisfies $\nabla_{\xi}\Theta(t,x,s,y,\xi)=0$ and from \eqref{nablaThet} and \eqref{determinhessphi} this translates into
\be
\Big((g^{j,k}(z_{x,y}))_{j,k}\Big) (\nabla_{x}\phi)^{t}(0,z_{x,y},\xi)=\frac{(x-y)}{(t-s)}.
\ee
Since $(g^{j,k})_{j,k}$ is invertible and using again \eqref{determinhessphi} we can apply the implicit function's theorem to obtain (for $T$ small enough) a critical point $\xi_{c}=\xi_{c}(t,x,s,y)$ for $\Theta$.
In order to show that $\xi_{c}$ is non-degenerate we compute
\begin{multline}
\partial_{\xi_{q}}\partial_{\xi_{l}}\Theta(t,x,s,y,\xi)=\sum_{j=1}^{n}\partial^{3}_{\xi_{q},\xi_{l},x_{j}}\phi(0,z_{x,y},\xi)\Big(2\sum_{k=1}^{n}g^{j,k}(z_{x,y})\partial_{x_{k}}\phi(0,z_{x,y},\xi)-\frac{(x_{j}-y_{j})}{(t-s)}\Big)\\
+2\sum_{j=1}^{n}\partial^{2}_{\xi_{l},x_{j}}\phi(0,z_{x,y},\xi)\Big(\sum_{k=1}^{n}g^{j,k}(z_{x,y})\partial^{2}_{\xi_{q},x_{k}}\phi(0,z_{x,y},\xi)\Big),
\end{multline}
consequently at the critical point $\xi=\xi_{c}$ the hessian matrix $\nabla^{2}_{\xi,\xi}\Theta$ is given by
\[
\nabla^{2}_{\xi,\xi}\Theta(t,x,s,y,\xi_{c})=2(\nabla^{2}_{\xi,x}\phi)(g^{ij}(z_{x,y}))_{i,j}(\nabla^{2}_{\xi,x}\phi)|_{(0,z_{x,y},\xi_{c})},
\]
and thereforee for $T$ small enough the critical point $\xi_{c}$ is non-degenerate by \eqref{determinhessphi}.
\end{proof}

On the support of $\kappa$ it follows that the kernel $K^{s}_{h}$ writes
\begin{multline}
K^{s}_{h}(t,x,s,y)=\frac{1}{(2\pi h)^{n}}\int e^{\frac{i}{h}|t-s|\Theta(t,x,s,y,\xi)}\psi^{2}(|\xi|)(1-\rho(h^{-1/3}|t-s|))\\ \times c_{m}(x,\xi/h)\overline{c_{m}(y,\xi/h)}d\xi,
\end{multline}
where, if $\omega=|t-s|/h$ and $\xi_{1}\simeq 1$, the symbol satisfies 
\[
|\partial^{\alpha}_{t,s}\partial^{k}_{\omega}\sigma_{h}(t,x,s,y,\omega\xi/|t-s|)|\leq C_{\alpha,k} h^{-|\alpha|/3}(|t-s|^{3/2}/h)^{-2k/3},
\]
where we set 
\[
\sigma_{h}(t,x,s,y,\omega\xi/|t-s|)=(1-\rho(h^{-1/3}|t-s|))c_{m}(x,\omega\xi/|t-s|)\overline{c_{m}(y,\omega\xi/|t-s|)}.
\]
Indeed, since $c_{m}\in S^{0}_{2/3,1/3}$, for $\alpha=0$ one has
\begin{eqnarray*}
|\partial^{k}_{\omega}\sigma_{h}| & \leq & |\xi||t-s|^{-k}|(\partial^{k}_{\xi}c_{m})(t,x,\omega\xi/|t-s|)|\\
 & \leq & C_{0,k}|t-s|^{-k}(\omega/|t-s|)^{-2k/3}\\
& = & C_{0,k}|t-s|^{-k}h^{2k/3}.
\end{eqnarray*}
We conclude using the next lemma with $\omega=\frac{|t-s|}{h}$ and $\delta=|t-s|^{3/2}\geq h^{1/2}\gg h$. 
\begin{lemma}\label{lemsta}
Suppose that $\Theta(z,\xi)\in C^{\infty}(\mathbb{R}^{2(n+1)}\times\mathbb{R}^{n})$ is real, $\nabla_{\xi}\Theta(z,\xi_{c}(z))=0$, $\nabla_{\xi}\Theta(z,\xi) \neq 0$ if $\xi\neq \xi_{c}(z)$, and 
\[
|\det \nabla^{2}_{\xi\xi}\Theta|\geq c_{0}>0, \ if\ |\xi|\leq 1. 
\] 
Suppose also that 
\[
|\partial^{\alpha}_{z}\partial^{\beta}_{\xi}\Theta(z,\xi)|\leq C_{\alpha,\beta} h^{-|\alpha|/3}, \quad \forall\alpha,\beta. 
\]
In addition, suppose that the symbol $\sigma_{h}(z,\xi,\omega)$ vanishes when $|\xi|\geq 1$ and satisfies 
\[
|\partial^{\alpha}_{z}\partial^{\gamma}_{\xi}\partial^{k}_{\omega}\sigma_{h}(z,\xi,\omega)|\leq C_{k,\alpha,\gamma} h^{-(|\alpha|+|\gamma|)/3}(\delta/h)^{-2k/3},\quad \forall k,\alpha,\gamma,
\]
where on the support of $\sigma_{h}$ we have $\omega\geq h^{-2/3}$ and $\delta>0$. Then we can write
\[
\int_{\mathbb{R}^{n}} e^{i\omega\Theta(z,\xi)}\sigma_{h}(z,\xi,\omega)d\xi=\omega^{-n/2}e^{i\omega\Theta(z,\xi_{c}(z))}b_{h}(z,\omega),
\]
where $b_{h}$ satisfies
\[
|\partial^{k}_{\omega}\partial^{\alpha}_{z}b_{h}(z,\omega)|\leq C_{k,\alpha} h^{-|\alpha|/3}(\delta/h)^{-2k/3}
\]
and where each of the constants depend only on $c_{0}$ and the size of finitely many of the constants $C_{\alpha,\beta}$ and $C_{k,\alpha,\gamma}$ above. In particular, the constants are uniform in $\delta$ if $1\geq \delta\geq h$. 
\end{lemma}
This Lemma is used in \cite[Lemma 2.6]{smso95} and also in the thesis of Grieser \cite{grie92} and it follows easily from the proof of the standard stationary phase lemma (see \cite[pag. 45]{so93}).
Proposition \ref{propfio} is thus proved.
\end{proof}
For each $t$, $s$, let $T^{s}_{h}(t,s)$ be the "frozen" operator defined by
\[
T^{s}_{h}(t,s)g(x)=\int K^{s}_{h}(t,x,s,y)g(y)dy.
\]
From Proposition \ref{propfio} we deduce
\be\label{linf}
\|T^{s}_{h}(t,s)g\|_{L^{\infty}(\mathbb{R}^{n})}\leq C\max(h^{-n},(h|t-s|)^{-n/2})\|g\|_{L^{1}(\mathbb{R}^{n})}.
\ee
We need the following 
\begin{lemma}
If $T$ is small enough then for $t$, $s$ fixed the frozen operators $T^{s}_{h}(t,s)$, $T^{f}_{h}(t,s)$ are bounded on $L^{2}(\mathbb{R}^{n})$, i.e. for all $g\in L^{2}(\mathbb{R}^{n})$ we have
\be\label{energie}
\|T^{s}_{h}(t,s)g\|_{L^{2}(\mathbb{R}^{n})}\leq C\|g\|_{L^{2}(\mathbb{R}^{n})}.
\ee
\end{lemma}
\begin{proof}
If $f\in L^{2}(\mathbb{R}^{n})$ then
\begin{multline}
\|W_{h}f(.,t)\|^{2}_{L^{2}(\mathbb{R}^{n})}=
\frac{1}{(2\pi h)^{2n}}\int_{\xi,\eta}\int_{x}e^{\frac{i}{h}(\phi(t,x,\xi)-\phi(t,x,\eta))}c_{m}(x,\xi/h)\overline{c_{m}(x,\eta/h)}\\ \times \psi(|\xi|)\psi(|\eta|)\hat{f}(\frac{\xi}{h})\hat{\bar{f}}(\frac{\eta}{h})dx d\xi d\eta.
\end{multline}
From Lemma \ref{lemdetermin} it folows that the mapping
\[
\chi:=\Big(x\rightarrow -t(\xi_{1}+\eta_{1},0,..,0)+\int_{0}^{1}\nabla_{\xi}\phi(0,x,(1-w)\xi+w\eta)dw\Big)
\]
is a diffeomorphisme away from the hypersurface $\zeta=0$ with uniform bounds of the Jacobian of $\chi^{-1}$. This change of variables reduces the problem to the $L^{2}$-continuity of semi-classical pseudo-differential operators with symbols of type $(2/3,1/3)$. 
\end{proof}
Interpolation between \eqref{linf} and \eqref{energie} with weights $1-2/r$ and $2/r$ respectively yields
\be
\|T^{s}_{h}(t,s)g\|_{L^{r}(\mathbb{R}^{n})}\leq Ch^{-n(1-2/r)}(1+\frac{|t-s|}{h})^{-n(1/2-1/r)}\|g\|_{L^{r'}(\mathbb{R}^{n})}
\ee
and hence
\[
\|T^{s}_{h}F\|_{L^{q}(0,T],L^{r}(\mathbb{R}^{n})}\leq C h^{-n/2(1-2/r)}\| \int_{1\ll \frac{|t-s|}{h}}^{T}|t-s|^{-n/2(1-2/r)}\|F(.,s)\|_{L^{r'}(\mathbb{R}^{n})}ds\|_{L^{q'}((0,T])}.
\]
Since $n(\frac{1}{2}-\frac{1}{r})=\frac{2}{q}<1$ the application $|t|^{-2/q}:L^{q'}\rightarrow L^{q}$ is bounded and by Hardy-Littlewood-Sobolev inequality we deduce 
\be
\|T^{s}_{h}F\|_{L^{q}((0,T],L^{r}(\mathbb{R}^{n}))}\leq Ch^{-\frac{2}{q}} \|F\|_{L^{q'}((0,T],L^{r'}(\mathbb{R}^{n}))}.
\ee

\item The "frozen" term $T^{f}_{h}$

To estimate $T^{f}_{h}$ it suffices to obtain bounds for its kernel $K^{f}_{h}$ with both the variables $(t,x)$ and $(s,y)$ restricted to lie in a cube of $\mathbb{R}^{n+1}$ of sidelength comparable to $h^{1/3}$. Let us decompose $S_{T}$ into disjoint cubes $Q=Q_{x}\times Q_{t}$ of sidelength $h^{1/3}$. We then have
\[
\|T^{f}_{h}F\|^{q}_{L^{q}([0,T],L^{r}(\mathbb{R}^{n}))}=\int_{0}^{T}\Big(\sum_{Q=Q_{x}\times Q_{t}}\|\chi_{Q}T^{f}_{h}F\|^{r}_{L^{r}(Q_{x})}\Big)^{q/r}dt=\sum_{Q}\|\chi_{Q}T^{f}_{h}F\|^{q}_{L^{q}([0,T],L^{r}(\mathbb{R}^{n}))},
\]
where by $\chi_{Q}$ we denoted the characteristic function of the cube $Q$.
In fact, by the definition, the integral kernel $K^{f}_{h}(t,x,s,y)$ of $T^{f}_{h}$ vanishes if $|t-s|\geq h^{1/3}$. If $|t-s|\leq h^{1/3}$ and $|x-y|\geq C_{0}h^{1/3}$, then the phase
\[
\phi(t,x,\xi)-\phi(s,y,\xi)
\]
has no critical points with respect to $\xi_{1}$ (on the support of $\psi$), so that
\[
|K^{f}_{h}(t,x,s,y)|\leq C_{N} h^{N}\quad \forall N,\ if\ |x-y|\geq C_{0} h^{1/3}.
\]
It thereforee suffices to estimate $\|\chi_{Q}T^{f}_{h}\chi_{Q^{*}}F\|_{L^{q}([0,T],L^{r}(\mathbb{R}^{n}))}$, where $Q^{*}$ is the dilate of $Q$ by some fixed factor independent of $h$. Since $q>2>q'$, $r\geq 2\geq r'$, where $q'$, $r'$ are such that $1/q+1/q'=1$, $1/r+1/r'=1$, then we shall obtain
\be\label{cub}
\sum_{Q}\|\chi_{Q}T^{f}_{h}\chi_{Q^{*}}F\|^{q}_{L^{q}([0,T],L^{r}(\mathbb{R}^{n}))}\leq C_{1}\sum_{Q}\|\chi_{Q^{*}}F\|^{q}_{L^{q'}([0,T],L^{r'}(\mathbb{R}^{n}))}\leq C_{2}\|F\|^{q}_{L^{q'}([0,T],L^{r'}(\mathbb{R}^{n}))}.
\ee
In order to prove \eqref{cub} we shall use the following:
\begin{prop}\label{propl2multip}
Let $b(\xi)\in L^{\infty}(\mathbb{R}^{n})$ be elliptic near $\xi_{1}\simeq 1$, $b_{h}(\xi):=b(\xi/h)$, then for $h\ll |t-s|\leq h^{1/3}$, $h\ll |x-y|\leq h^{1/3}$ the operator defined by
\be
B_{h}f(x,t)=\frac{1}{(2\pi h)^{n}}\int e^{\frac{i}{h}\phi(t,x,\xi)}\psi(|\xi|)b_{h}(\xi)\hat{f}(\frac{\xi}{h})d\xi
\ee
satisfies
\be
\|B_{h}f\|_{L^{q}((0,T],L^{r}(\mathbb{R}^{n}))}\leq C h^{-\frac{1}{q}}\|f\|_{L^{2}(\mathbb{R}^{n})}.
\ee
\end{prop}
\begin{proof}
We use again the $\text{TT}^{*}$ argument. Since $b(\xi)$ acts as an $L^{2}$ multiplier we can apply the stationary phase theorem in the integral
\[
\int e^{\frac{i}{h}(\phi(t,x,\xi)-\phi(s,y,\xi))}\psi(|\xi|)d\xi
\]
in order to obtain
\[
\|B_{h}B^{*}_{h}F\|_{L^{q}((0,T],L^{r}(\mathbb(\mathbb{R}^{n}))}\lesssim h^{-\frac{2}{q}}\|F\|_{L^{q'}((0,T],L^{r'}(\mathbb{R}^{n}))}.
\]
Notice that we haven't used the special properties of the phase function at $t=0$.
\end{proof}
Let now $Q$ be a fixed cube in $\mathbb{R}^{n+1}$ of sidelength $h^{1/3}$. Let
\[
b_{h}(t,x,s,y,\xi)=\rho(h^{-1/3}|t-s|)c_{m}(x,\xi/h)\overline{c_{m}(y,\xi/h)},
\]
and write
\begin{multline}
  b_{h}(t,x,s,y,\xi)=b_{h}(0,0,s,y,\xi)+\int_{0}^{t}\partial_{t}b_{h}(r,0,s,y,\xi)dr\\
{}+\cdots +\int_{0}^{t}\cdots\int_{0}^{x_{n}}\partial_{t}\cdots\partial_{x_{n}}b_{h}(r,z_{1},..,z_{n},s,y,\xi)drdz.
\end{multline}
If the symbol $c$ is independent of $t$, $x$ then the estimates \eqref{eststri}
follow from Proposition \ref{propl2multip}. We use this, for instance, to deduce 
\begin{multline}
  \|\chi_{Q}T^{f}_{h}\chi_{Q*}F\|_{L^{q}((0,T],L^{r}(\mathbb{R}^{n}))}\leq
  Ch^{-n/2(1/2-1/r)}\\
\times \left(\left\|\int\int e^{\frac{i}{h}(x\xi-\phi(s,y,\xi))}\psi(|\xi|)b_{h}(0,0,s,y,\xi)F(y,s)d\xi ds dy\right\|_{L^{2}(\mathbb{R}^n)}\right.\\
\left.{}+\cdots+\int_{0}^{h^{1/3}}\int_{0}^{h^{1/3}}\left\|\int\int e^{\frac{i}{h}(x\xi-\phi(s,y,\xi))}\partial_{t}..\partial_{x_{n}}\psi(|\xi|)b_{h}(r,z,s,y,\xi)F(y,s)d\xi ds dy\right\|_{L^{2}(\mathbb{R}^{n})}drdz\right).
\end{multline}
Each derivative of $b_{h}(t,x,s,y,\xi)$ loses a factor of $h^{-1/3}$, but this is compensated by the integral over $(r,z)$, so that it suffices to establish uniform estimates for fixed $(r,z)$. By duality, we have to establish the estimate
\[
\|\int\int e^{\frac{i}{h}\phi(s,y,\xi)}\psi(|\xi|)b_{h}(0,0,s,y,\xi)\hat{f}(\frac{\xi}{h})d\xi\|_{L^{q}((0,T],L^{r}(\mathbb{R}^{n}))}\leq C\|f\|_{L^{2}(\mathbb{R}^{n})},
\] 
which follows by using the same argument of freezing the variables $(s,y)$ together with the Proposition \ref{propl2multip}.
\end{itemize}

\subsubsection{The diffractive term $D_{h}$}\label{sectdiffterm}
In order to estimate the diffractive term we shall proceed again like in \cite[Sect.2]{smso95}. 
\begin{lemma}\label{lemestim}
For $x_{n}\geq 0$ and for $\xi$ in a small conic neighborhood of the positive $\xi_{1}$ axis, the symbol $q$ of $S_{h}$ can be written in the form
\begin{eqnarray*}
q(x,\xi) & := &(a(x,\xi)((1-\chi)A_{+})(\zeta(x,\xi))+b(x,\xi)((1-\chi)A_{+})'(\zeta(x,\xi)))\frac{Ai(\zeta_{0}(\xi))}{A_{+}(\zeta_{0}(\xi))}\\
& = & p(x,\xi,\zeta(x,\xi)),
\end{eqnarray*}
where, for some $c>0$
\[
|\partial^{\alpha}_{\xi}\partial^{j}_{\zeta}\partial^{\beta}_{x'}\partial^{k}_{x_{n}} p(x,\xi,\zeta(x,\xi))|\leq C_{\alpha,j,\beta,k}\xi^{1/6-|\alpha|+2k/3}_{1} e^{-c x_{n}^{3/2}\xi_{1}-|\zeta|^{3/2}/2}.
\]
\end{lemma}
\begin{proof}
Since 
\[
|\partial^{k}_{\zeta}((1-\chi)A_{+})(\zeta)|\leq C_{k,\epsilon} e^{(2/3+\epsilon)|\zeta|^{3/2}},\quad \forall\epsilon>0,
\]
and the symbols $a$ and $b$ belong to $S^{1/6}_{1,0}$, the above fact will follow by showing that in the region $\zeta(x,\xi)\geq -2$,
\[
\frac{Ai}{A_{+}}(\zeta_{0}(\xi))=\tilde{p}(x,\xi',\zeta(x,\xi)),
\]
where if $\xi'=(\xi_{1},..,\xi_{n-1})$
\be\label{diff}
|\partial^{\alpha}_{\xi'}\partial^{j}_{\zeta}\partial^{\beta}_{x'}\partial^{k}_{x_{n}} \tilde{p}(x,\xi',\zeta)|\leq C_{\alpha,j,\beta,k,\epsilon}\xi^{-|\alpha|+2k/3}_{1}e^{-cx^{3/2}_{n}\xi_{1}-(4/3-\epsilon)|\zeta|^{3/2}}.
\ee
At $x_{n}=0$, one has $\zeta=\zeta_{0}$, $\partial_{x_{n}}\zeta<0$. It follows that for some $c>0$
\[
\zeta_{0}(x,\xi)\geq \zeta(x,\xi)+cx_{n}\xi^{2/3}_{1}.
\]
By the asymptotic behavior of the Airy function we have, in the region $\zeta(x,\xi)\geq -2$
\be\label{estai}
|\Big(\frac{Ai}{A_{+}}\Big)^{(k)}(\zeta_{0})|\leq C_{k,\epsilon} e^{-c x^{3/2}_{n}\xi_{1}-(4/3-\epsilon)|\zeta(x,\xi)|^{3/2}}.
\ee
We introduce a new variable $\tau(x,\xi)=\xi^{1/3}_{1}\zeta(x,\xi)$.
At $x_{n}=0$ one has $\tau=-\xi_{n}$, so that we can write $\xi_{n}=\sigma(x,\xi',\tau)$, where $\sigma$ is homogeneous of degree $1$ in $(\xi',\tau)$. We set
\[
\tilde{p}(x,\xi',\zeta)=\frac{Ai}{A_{+}}(-\xi^{-1/3}_{1}\sigma(x,\xi',\xi^{1/3}\zeta)).
\]
The estimates \eqref{diff} will follow by showing that
\be\label{estdiff}
|\partial^{\alpha}_{\xi'}\partial^{j}_{\tau}\partial^{\beta}_{x'}\partial^{k}_{x_{n}}\frac{Ai}{A_{+}}(-\xi^{-1/3}_{1}\sigma(x,\xi',\tau))|\leq C_{\alpha,j,\beta,k,\epsilon}\xi^{-|\alpha|-j+2k/3}_{1} e^{-cx^{3/2}_{n}\xi_{1}-(4/3-\epsilon)|\tau|^{3/2}\xi^{-1/2}_{1}}.
\ee
For $k=0$, the estimates \eqref{estdiff} follow from \eqref{estai}, together with the fact that
\[
|\partial^{\alpha}_{\xi'}\partial^{j}_{\tau}\partial^{\beta}_{x'}\frac{Ai}{A_{+}}(-\xi^{-1/3}_{1}\sigma(x,\xi',\tau))|\leq C_{\alpha,\beta,j}(x_{n}\xi^{2/3}_{1}+\xi^{-1/3}_{1}|\tau|)\xi^{-|\alpha|-j}_{1},
\]
which, in turn, holds by homogeneity, together with the fact that $\sigma(x,\xi',\tau)=0$ if $x_{n}=\tau=0$. If $k>0$, the estimate \eqref{estdiff} follows by observing that the effect of differentiating in $x_{n}$ is similar to multiplying by a symbol of order $2/3$. This concludes the proof of Lemma \ref{lemestim}.
\end{proof}
\begin{lemma}
The Schwartz kernel of the diffractive term $D_{h}$ writes in the form
\begin{multline}
  \label{formpha0}
  \int e^{i(\theta(x,\xi)-ht\xi^{2}_{1})}\psi(h|\xi|)q(x,\xi)d\xi\\
=\int e^{i(\theta(x,\xi)-ht\xi^{2}_{1}+\sigma \xi^{-2/3}_{1}\zeta(x,\xi)+\sigma^{3}/3\xi^{2}_{1}-<y,\xi>)}\psi(h|\xi|)c_{d}(x,\xi,\sigma)d\sigma d\xi,
\end{multline}
where $<.,.>$ denotes the scalar product and where
\[
|\partial^{\alpha}_{\xi}\partial^{j}_{\sigma}\partial^{\beta}_{x'}\partial^{k}_{x_{n}}c_{d}(x,\xi,\sigma)|\leq C_{\alpha,j,\beta,k,N}\xi_{1}^{-1/2-|\alpha|-2j/3+2k/3}e^{-cx_{n}^{3/2}\xi_{1}}(1+\xi^{-4/3}_{1}\sigma^{2})^{-N/2},\quad\forall N.
\]
\end{lemma}
\begin{proof}
The symbol $c_{d}$ of the Schwartz kernel of $D_{h}$ writes as a product of two symbols $$c_{d}(x,\xi,\sigma)=c_{1}(x,\xi,\sigma\xi_{1}^{-2/3})c_{2}(x,\xi,\zeta(x,\xi)),$$ where
\[
c_{1}(x,\xi,\sigma\xi_{1}^{-2/3})=\xi_{1}^{-2/3}\Psi_{+}(\xi_{1}^{-2/3}\sigma)(a(x,\xi)+\sigma\xi_{1}^{-2/3}b(x,\xi))\in S^{-1/2}_{2/3,1/3}(\mathbb{R}^{n}_{x},\mathbb{R}^{n+1}_{\xi,\sigma})
\]  
comes from the Fourier transform of $A_{+}$ (here $\Psi_{+}$ is a symbol of order $0$) and where $c_{2}$ satisfies for all $N\geq 0$ (for $\sigma^{2}\xi_{1}^{-4/3}+\zeta(x,\xi)=0$)
\be\label{simbc2}
|\partial^{\alpha}_{\xi'}\partial^{j}_{\sigma}\partial^{\beta}_{x'}\partial^{k}_{x_{n}}c_{2}(x,\xi',-(\sigma^{2}\xi_{1}^{-4/3}))|\leq C_{\alpha,j,\beta,k,N}\xi_{1}^{-2j/3}|\sigma\xi_{1}^{-2/3}|^{j}\xi_{1}^{-|\alpha|+2k/3}e^{-cx_{n}^{3/2}\xi_{1}}(1+\xi_{1}^{-4/3}\sigma^{2})^{-N/2},
\ee
which follows from \eqref{diff}. We use the exponential factor $e^{-cx_{n}^{3/2}\xi_{1}}$ to deduce from \eqref{simbc2}
\[
\forall N,\quad |x_{n}^{j}\partial^{k}_{x_{n}}c_{2}(x,\xi',-(\sigma^{2}\xi_{1}^{-4/3}))|\leq C_{j,k,N}(x_{n}\xi_{1}^{2/3})^{j}e^{-c(x_{n}\xi_{1}^{2/3})^{3/2}}\xi_{1}^{2/3(k-j)}(1+\xi_{1}^{-4/3}\sigma^{2})^{-N/2},
\]
\end{proof}
From now on we proceed as for the main term and we reduce the problem to considering the operator
\[
W^{d}_{h}f(x,t)=\frac{1}{(2\pi h)^{n}}\int e^{\frac{i}{h}\tilde{\phi}(t,x,\xi,\sigma)}c_{d}(x,\xi/h,\sigma)\psi(|\xi|)\hat{f}(\frac{\xi}{h})d\xi,
\]
where $x_{n}^{j}\partial^{k}_{x_{n}}c_{d}\in S^{2(k-j)/3}_{2/3,1/3}(\mathbb{R}^{n-1}_{x'}\times\mathbb{R}^{n}_{\xi})$ uniformly over $x_{n}$ and where we set
\begin{equation}\label{phasatilde}
\tilde{\phi}(t,x,\xi,\sigma):=-t\xi^{2}_{1}+\theta(x,\xi)+\sigma\xi^{1/3}_{1}\zeta(x,\xi)+\xi_{1}\sigma^{3}/3,
\end{equation}
obtained after the changes of variables $\sigma\rightarrow\sigma\xi_{1}$, $\xi\rightarrow\xi/h$ in  \eqref{formpha0}. Using the freezing arguments behind the proof of the estimates for $T^{f}_{h}$ and Minkovski inequality we have
\begin{multline*}
\|W^{d}_{h}f\|_{L^{q}((0,T],L^{r}(\mathbb{R}^{n}))}\leq\left\|\frac{1}{(2\pi h)^{n}}\int e^{\frac{i}{h}\tilde{\phi}(t,x,\xi,\sigma)}c_{d}(x',0,\xi/h,\sigma)\psi(|\xi|)\hat{f}(\frac{\xi}{h})d\sigma d\xi\right\|_{L^{q}((0,T],L^{r}(\mathbb{R}^{n}))}\\
{}+h^{-2/3}\int_{0}^{h^{2/3}}\left\|\frac{1}{(2\pi h)^{n}}\int e^{\frac{i}{h}\tilde{\phi}(t,x,\xi,\sigma)}h^{2/3}\partial_{x_{n}}c_{d}(x',r,\xi/h,\sigma)\psi(|\xi|)\hat{f}(\frac{\xi}{h})d\sigma d\xi\right\|_{L^{q}((0,T],L^{r}(\mathbb{R}^{n-1}))}dr\\
{}+h^{2/3}\int_{r>h^{2/3}}\frac{dr}{r^{2}}\left\|\frac{1}{(2\pi h)^{n}}\int e^{\frac{i}{h}\tilde{\phi}(t,x,\xi,\sigma)}h^{-2/3}r^{2}\partial_{x_{n}}c_{d}(x',r,\xi/h,\sigma)\psi(|\xi|)\hat{f}(\frac{\xi}{h})d\sigma d\xi\right\|_{L^{q}((0,T],L^{r}(\mathbb{R}^{n-1}))}.
\end{multline*}
Since $c_{d}(x',0,\xi,\sigma)$ and $h^{2/3}(1+h^{-4/3}r^{2})\partial_{x_{n}}c_{d}(x',r,\xi,\sigma)$ are symbols of order $0$ and type $(2/3,1/3)$ with uniform estimates over $r$, the estimates for the diffractive term also follow from Proposition \ref{propestim}. Indeed, the term in the second line loses a factor $h^{-2/3}$, but this is compensated by the integral over $r\leq h^{2/3}$. The term in the third line can be bounded by above by 
\begin{multline*}
 h^{2/3}\int_{r>h^{2/3}}\frac{dr}{r^{2}}\times\left\|\frac{1}{(2\pi h)^{n}}\int e^{\frac{i}{h}\tilde{\phi}(t,x,\xi,\sigma)}(h^{-2/3}r^{2}\partial_{x_{n}}c_{d}(x',r,\xi/h,\sigma))\psi(|\xi|)\hat{f}(\frac{\xi}{h})d\sigma d\xi\right\|_{L^{q}((0,T],L^{r}(\mathbb{R}^{n}))}\\
\leq\left\|\frac{1}{(2\pi h)^{n}}\int e^{\frac{i}{h}\tilde{\phi}(t,x,\xi,\sigma)}(h^{-2/3}r^{2}\partial_{x_{n}}c_{d}(x',r,\xi/h))\psi(|\xi|)\hat{f}(\frac{\xi}{h})d\sigma d\xi\right\|_{L^{q}((0,T],L^{r}(\mathbb{R}^{n}))}.
\end{multline*}
We conclude using the same arguments as in the proof of Proposition \ref{propestim}, where now $W_{h}$ is replaced by operators with symbols $c_{d}(x',0,\xi,\sigma)$, $h^{2/3}\partial_{x_{n}}c_{d}(x',r,\xi,\sigma)$ and $h^{-2/3}r^{2}\partial_{x_{n}}c_{d}(x',r,\xi,\sigma)$ respectively.
Notice however that for this term we can't apply directly Lemma \ref{lemptcrit} since the expansion of the Airy function giving the phase function \eqref{phi} is available only for $\zeta(x,\xi/h)\leq -1$. Writing the phase function of \eqref{formpha0} in the form $\tilde{\phi}(t,x,\xi,\sigma)-<y,\xi>$,
we notice that at $t=0$ this phase is homogeneous of degree $1$ in $\xi$ and the proof of the non-degeneracy of the critical points in the $\text{TT}^{*}$ argument of Lemma \ref{lemptcrit} reduces to checking that the Jacobian $J$ of the mapping
\be\label{jacaplic}
(\xi,\sigma)\rightarrow (\nabla_{x}(\theta(x,\xi)+\sigma\zeta(x,\xi)),\zeta(x,\xi)+\sigma^{2})
\ee
does not vanish at the critical point of the phase of \eqref{formpha0}. Hence we will obtain a phase function $\breve{\phi}(t,x,\xi)$ which will satisfy $\nabla^{2}_{x,\xi}\breve{\phi}(0,x,\xi)\neq 0$ and this will hold also for small $|t|\leq T$ and we can use the same argument as in Lemma \ref{lemptcrit}.
To prove that the Jacobian of the application \eqref{jacaplic} doesn't vanish we use \cite[Lemma A.2]{smso94}. Precisely, at this (critical) point $\sigma=\zeta(x,\xi)=0$, $y=0$ and $\nabla_{x'}\zeta(x,\xi)=0$. Since $\partial_{x_{n}}\zeta(x,\xi)\neq 0$ and $\partial_{\xi_{n}}\zeta(x,\xi)\neq 0$ there, the result follows by the nonvanishing of $|\nabla_{x'}\nabla_{\xi'}\theta(x,\xi)|$. In fact we have
\[
\det
\left(
\begin{array}{ccc}
   \nabla_{x'}\nabla_{\xi'}\theta&  \nabla_{\xi'}\partial_{x_{n}}\theta& \nabla_{\xi'}\zeta  \\
  \partial_{\xi_{n}}\nabla_{x'}\theta& \partial_{\xi_{n}}\partial_{x_{n}}\theta   & \partial_{\xi_{n}}\zeta   \\
  \nabla_{x'}\zeta& \partial_{x_{n}}\zeta & 2\sigma
\end{array}
\right)|_{\sigma^{2}=-\zeta=0}\neq 0.
\]

\section{Strichartz estimates for the classical Schr\"odinger \\ equation outside a strictly convex obstacle in $\mathbb{R}^{n}$}\label{secclasic}
In this section we prove Theorem \ref{thmstri} under the Assumptions \ref{assum2}. In what follows we shall work with the Laplace operator with constant coefficients $\Delta_{D}=\sum_{j=1}^{n}\partial^{2}_{j}$ acting on $L^{2}(\Omega)$ to avoid technicalities, where $\Omega$ is the exterior in $\mathbb{R}^{n}$ of a strictly convex domain $\Theta$. In the proof of Theorem \ref{thmstri} we distinguish two main steps: we start by performing a time rescaling which transforms the equation \eqref{schrod} into a semi-classical problem: due to the finite speed of propagation (proved by Lebeau \cite{le92}), we can use the (local) semi-classical result of Theorem \ref{thmstrichartz} together with the smoothing effect (following Staffilani and Tataru \cite{stta02} and Burq \cite{bu02}) to obtain classical Strichartz estimates near the boundary. Outside a fixed neighborhood of $\partial\Omega$ we use a method suggested by Staffilani and Tataru \cite{stta02} which consists in considering the Schr\"odinger flow as a solution of a problem in the whole space $\mathbb{R}^{n}$, for which the Strichartz estimates are known.

We start by proving that using Theorem \ref{thmstrichartz} on a compact manifold with strictly concave boundary we can deduce sharp Strichartz estimates for the semi-classical Schr\"odinger flow on $\Omega$. Precisely, the following holds
\begin{prop}\label{propstri}
Given $(q,r)$ satisfying the scaling condition \eqref{scaling} with $q>2$ there exists a constant $C>0$ such that the (classical) Schr\"odinger flow on $\Omega\times \mathbb{R}$ with Dirichlet boundary condition and spectrally localized initial data $\Psi(-h^{2}\Delta_{D})u_{0}$, where $\Psi\in C^{\infty}_{0}(\mathbb{R}\setminus 0)$, satisfies  
\be\label{estvo}
\|e^{it\Delta_{D}}\Psi(-h^{2}\Delta_{D})u_{0}\|_{L^{q}(\mathbb{R})L^{r}(\Omega)}\leq 
C\|\Psi(-h^{2}\Delta_{D})u_{0}\|_{L^{2}(\Omega)}.
\ee
\end{prop}
\begin{rmq}
We first proceed with the proof of Proposition \ref{propstri} and then we show how it can be used to prove Theorem \ref{thmstri}. For the proof of Proposition \ref{propstri} we use a similar method as the one given in our recent paper \cite{ip09} in collaboration with F.Planchon.
\end{rmq}

\begin{proof}
Let $\tilde{\Psi}\in C^{\infty}_{0}(\mathbb{R}\setminus \{0\})$ be
such that $\tilde{\Psi}=1$ on the support of $\Psi$, hence 
$$
\tilde{\Psi}(-h^{2}\Delta_{D})\Psi(-h^{2}\Delta_{D})=\Psi(-h^{2}\Delta_{D}).
$$
Following \cite{bu02}, \cite{ip09}, we split $e^{it\Delta_{D}}\Psi(-h^{2}\Delta_{D})u_{0}(x)$ as a sum of two terms 
\[
\tilde{\Psi}(-h^{2}\Delta_{D})\chi \Psi(-h^{2}\Delta_{D}) e^{it\Delta_{D}}u_{0}+\tilde{\Psi}(-h^{2}\Delta_{D})(1-\chi)\Psi(-h^{2}\Delta_{D})e^{it\Delta_{D}}u_{0},
\]
where $\chi\in C^{\infty}_{0}(\mathbb{R}^{n})$ equals $1$ in a neighborhood of $\partial\Omega$.

\begin{itemize}
\item Study of ${\tilde\Psi(-h^{2}\Delta_{D})}(1-\chi)\Psi(-h^{2}\Delta_{D})  e^{it\Delta_{D}}u_{0}$ :
    
Set $w_{h}(x,t)=(1-\chi)\Psi(-h^{2}\Delta_{D}) e^{it\Delta_{D}}u_{0}(x)$. Then $w_{h}$ satisfies
\begin{equation}\label{estloin1}
\left\{\begin{array}{ll}
i\partial_{t}w_{h}+\Delta_{D}w_{h}=-[\Delta_{D},\chi]\Psi(-h^{2}\Delta_{D}) e^{it\Delta_{D}}u_{0},\\
w_{h}|_{t=0}=(1-\chi)\Psi(-h^{2}\Delta_{D}) u_{0}.
\end{array}
\right.
\end{equation}
Since $\chi$ is equal to $1$ near the boundary $\partial\Omega$, the solution to \eqref{estloin1} solves also a
problem in the whole space $\mathbb{R}^{n}$. Consequently, the Duhamel
formula writes
\begin{equation}\label{est1}
w_{h}(t,x)=e^{it\Delta}(1-\chi)\Psi(-h^{2}\Delta_{D}) u_{0}-\int_{0}^{t} e^{i(t-s)\Delta}[\Delta_{D},\chi]\Psi(-h^{2}\Delta_{D}) e^{it\Delta_{D}}u_{0}(s)ds,
\end{equation}
where by $\Delta$ we denoted the free Laplacian on $\mathbb{R}^{n}$ and thereforee
the contribution of $e^{it\Delta}(1-\chi)\Psi(-h^{2}\Delta_{D}) u_{0}$
satisfies the usual Strichartz estimates. For the second term in the right hand side of \eqref{est1} we use the next lemma, due to Christ and Kiselev \cite{chki01}:
\begin{lemma}(Christ and Kiselev)
Consider a bounded operator 
\[
T:L^{q'}(\mathbb{R},B_{1})\rightarrow L^{q}(\mathbb{R}, B_{2})
\]
given by a locally integrable kernel $K(t,s)$ with values in bounded operators from $B_{1}$ to $B_{2}$, where $B_{1}$ and $B_{2}$ are Banach spaces. Suppose that $q'<q$. Then the operator
\[
\tilde{T}f(t)=\int_{s<t} K(t,s)f(s)ds
\]
is bounded from $L^{q'}(\mathbb{R},B_{1})$ to $L^{q}(\mathbb{R},B_{2})$ and
\[
\|\tilde{T}\|_{L^{q'}(\mathbb{R},B_{1})\rightarrow L^{q}(\mathbb{R},B_{2})}\leq C (1-2^{-(1/q-1/q')})^{-1}\|T\|_{L^{q'}(\mathbb{R},B_{1})\rightarrow L^{q}(\mathbb{R},B_{2})}.
\]
\end{lemma} 
This lemma allows (since $q>2$) to replace the study of the second term in the right hand side of \eqref{est1} by that of
\[
\int_{0}^{\infty} e^{i(t-s)\Delta}[\Delta_{D},\chi]\Psi(-h^{2}\Delta_{D}) e^{it\Delta_{D}}u_{0}(s)(s)ds=:U_{0}U^{*}_{0}f(x,t),
\] 
where $U_{0}= e^{it\Delta}$ is bounded from $L^{2}(\mathbb{R}^{n})$ to $L^{q}(\mathbb{R},L^{r}(\mathbb{R}^{n}))$ and $U^{*}_{0}$ is bounded from $L^{2}(\mathbb{R},H^{-1/2}_{\text{comp}})$ to $L^{2}(\mathbb{R}^{n})$ and where we set
$f:=[\Delta_{D},\chi]\Psi(-h^{2}\Delta_{D}) e^{it\Delta_{D}}u_{0}$ 
which belongs to $L^{2}H^{-1/2}_{\text{comp}}(\Omega)$ by \cite[Prop.2.7]{bgt03}. 
The estimates for $w_{h}$ follow like in \cite{bgt03} and we find
\begin{multline}\label{whestqr}
\|w_{h}\|_{L^{q}(\mathbb{R},L^{r}(\Omega))}\leq C\|(1-\chi)\Psi(-h^{2}\Delta_{D}) u_{0}\|_{L^{2}(\mathbb{R}^{n})}\\
+\|[\Delta_{D},\chi]\Psi(-h^{2}\Delta_{D}) e^{it\Delta_{D}}u_{0}\|_{L^{2}(\mathbb{R},H^{-1/2}_{\text{comp}}(\Omega))}.
\end{multline}
The last term in \eqref{whestqr} can be estimated using \cite[Prop.2.7]{bgt03} by 
\be
C\|\Psi(-h^{2}\Delta_{D}) e^{it\Delta_{D}}u_{0}\|_{L^{2}(\mathbb{R},H^{1/2}_{\text{comp}}(\Omega))}\leq C\|\Psi(-h^{2}\Delta_{D})u_{0}\|_{L^{2}(\Omega)}.
\ee
Finally, we conclude this part using \cite[Thm.1.1]{ivpl08} which gives
\be
\|\Psi(-h^{2}\Delta_{D}) w_{h}\|_{L^{r}(\Omega)}\leq \|w_{h}\|_{L^{r}(\Omega)}.
\ee

\item Study of $\tilde{\Psi}(-h^{2}\Delta_{D})\chi\Psi(-h^{2}\Delta_{D})e^{it\Delta_{D}}u_{0}$:
 
Let $\varphi\in C^{\infty}_{0}((-1,2))$ equal to $1$ on $[0,1]$. For
$l\in\mathbb{Z}$ let 
\begin{equation}
  v_{h,l}=\varphi(t/h-l)\chi\Psi(-h^{2}\Delta_{D})e^{it\Delta_{D}}u_{0},
\end{equation}
 which is a solution to
\be\label{eqvh}
\left\{\begin{array}{ll}
i\partial_{t}v_{h,l}+\Delta_{D}v_{h,l}=\Big(\varphi(t/h-l)[\Delta_{D},\chi]+i\frac{\varphi'(t/h-l)}{h}\chi\Big)\Psi(-h^{2}\Delta_{D})e^{it\Delta_{D}}u_{0},\\
v_{h,l}|_{t<hl-h}=0,\quad v_{h,l}|_{t>hl+2h}=0.
\end{array}
\right.
\ee
We denote by $V_{h,l}$ the right-hand side of \eqref{eqvh}, so that
\be\label{defVhl1}
V_{h,l}=\Big(\varphi(t/h-l)[\Delta_{D},\chi]+i\frac{\varphi'(t/h-l)}{h}\chi\Big)\Psi(-h^{2}\Delta_{D})e^{it\Delta_{D}}u_{0}.
\ee
Let $Q\subset\mathbb{R}^{n}$ be an open cube sufficiently large such
that $\partial\Omega$ is contained in the interior of $Q$. We denote
by $S$ the punctured torus obtained from removing the obstacle $\Theta$ (recall that
$\Omega=\mathbb{R}^{n}\setminus\Theta$) in the compact manifold
obtained from $Q$ with periodic boundary conditions on $\partial
Q$. Notice that defined in this way $S$ coincides with the Sina\"i
billiard. Let $\Delta_{S}:=\sum_{j=1}^{n}\partial^{2}_{j}$ denote
the Laplace operator on the compact domain $S$.

On $S$, we may define a spectral localization operator using
eigenvalues $\lambda_k$ and eigenvectors $e_k$ of $\Delta_S$: if $f=\sum_k c_k e_k$, then
\begin{equation}
  \label{eq:DeltaS}
  \Psi(-h^{2}\Delta_{S})f=\sum_k \Psi(-h^{2}\lambda_k^2)
  c_k e_k.
\end{equation}
\begin{rmq}\label{rmqeqsd}
Notice that in a neighborhood of the boundary, the domains of $\Delta_{S}$ and $\Delta_{D}$ coincide, thus if $\tilde{\chi}\in C^{\infty}_{0}(\mathbb{R}^{n})$ is supported near $\partial\Omega$ then 
$\Delta_{S}\tilde{\chi}=\Delta_{D}\tilde{\chi}$. 
\end{rmq}

In what follows let $\tilde{\chi}\in C^{\infty}_{0}(\mathbb{R}^{n})$
be equal to $1$ on the support of $\chi$ and be supported in a
neighborhood of $\partial\Omega$ such that on its support the operator
$-\Delta_{D}$ coincide with  $-\Delta_{S}$. From their respective
definition, $v_{h,l}=\tilde{\chi}v_{h,l}$, $V_{h,l}=\tilde{\chi}V_{h,l}$,
consequently $v_{h,l}$ will also solve the following equation on the compact domain $S$ 
\begin{equation}\label{eqvhs}
\left\{\begin{array}{ll}
i\partial_{t}v_{h,l}+\Delta_{S}v_{h,l}=V_{h,l},\\
v_{h,l}|_{t<h(l-1/2)\pi}=0,\quad v_{h,l}|_{t>h(l+1)\pi}=0.
\end{array}
\right.
\end{equation}
Writing the Duhamel formula for the last equation \eqref{eqvhs} on $S$, applying $\tilde{\Psi}(-h^{2}\Delta_{D})$ and using that $\tilde{\chi}v_{h,l}=v_{h,l}$, $\tilde{\chi}V_{h,l}=V_{h,l}$ and writing
\begin{multline}
\tilde{\Psi}(-h^{2}\Delta_{D})\tilde{\chi}=\chi_{1}\tilde{\Psi}(-h^{2}\Delta_{S})\tilde{\chi}+(1-\chi_{1})\tilde{\Psi}(-h^{2}\Delta_{D})\tilde{\chi}\\+\chi_{1}(\tilde{\Psi}(-h^{2}\Delta_{D})-\tilde{\Psi}(-h^{2}\Delta_{S}))\tilde{\chi}
\end{multline}
for some $\chi_{1}\in C^{\infty}_{0}(\mathbb{R}^{n})$ equal to one on the support of $\tilde{\chi}$, yields
\begin{multline}\label{egalsd}
\tilde{\Psi}(-h^{2}\Delta_{D})v_{h,l}(x,t)
=\chi_{1}\int_{hl-l}^{t} e^{i(t-s)\Delta_{S}}\tilde{\Psi}(-h^{2}\Delta_{S})V_{h,l}(x,s)ds+\\
+(1-\chi_{1})\int_{hl-l}^{t} \tilde{\Psi}(-h^{2}\Delta_{D}) e^{i(t-s)\Delta_{S}}V_{h,l}(x,s)ds\\
+\chi_{1}(\tilde{\Psi}(-h^{2}\Delta_{D})-\tilde{\Psi}(-h^{2}\Delta_{S}))v_{h,l}.
\end{multline}
Denote by $v_{h,l,m}$ the first term of \eqref{egalsd}, by $v_{h,l,f}$
the second one and by $v_{h,l,s}$ the last one. We deal we them
separately. To estimate the $L^{q}_{t}L^{r}(\Omega)$ norm of
$v_{h,l,f}$ we notice that it is supported away from the boundary,
therefore the estimates will follow as in the previous part of this
section. Indeed, notice that since $v_{h,l}$ solves also the equation
\eqref{eqvh} on $\Omega$ we can use the Duhamel formula on $\Omega$ so
that in the integral defining $v_{h,l,f}$ to have $\Delta_{D}$ instead
of $\Delta_{S}$. We then estimate the $L^{q}_{t}L^{r}(\Omega)$ norm of
$v_{h,l,f}$ applying the Minkovski inequality and using the sharp
Strichartz estimates for $(1-\chi_{1})\tilde{\Psi}(-h^{2}\Delta_{D})
e^{i(t-s)\Delta_{D}}V_{h,l}$ deduced in the first part of the proof of
Proposition \ref{propstri} and obtain, denoting $I^h_l=[hl-h,hl+2h]$,
\be\label{estvhleeee}
\|v_{h,l,f}\|_{L^{q}(I^h_l ,L^{r}(\Omega))}\leq C\int_{I^h_l }\|V_{h,l}(x,s)\|_{L^{2}(\Omega)}ds.
\ee
For the last term $v_{h,l,s}$ we use the next lemma that will be proved in the Appendix:
\begin{lemma}\label{lemschwartzrest} Let $\chi_{1}\in C^{\infty}_{0}(\mathbb{R}^{n})$ be equal to $1$ on a fixed neighborhood of the support of $\tilde{\chi}$. Then we have
\be\label{estvhleeees}
\|v_{h,l,s}\|_{L^{q}(I^h_l ,L^{r}(\Omega))}\leq C_{N}h^{N}\|V_{h,l}(x,s)\|_{L^{2}(I^h_l ,H^{n(\frac{1}{2}-\frac{1}{r})-\frac{1}{2}}_{0}(\Omega))},\quad \forall N\in\mathbb{N}.
\ee
\end{lemma}
To estimate the main contribution $v_{h,l,m}$ we use the Minkovski inequality which yields
\begin{multline}\label{vhlnorm}
\|v_{h,l,m}\|_{L^{q}(I^h_l ,L^{r}(\Omega))}=\|v_{h,l,m}\|_{L^{q}(I^h_l ,L^{r}(S))}\\ \leq  C\int_{I^h_l }\|e^{i(t-s)\Delta_{S}}\tilde{\Psi}(-h^{2}\Delta_{S})V_{h,l}(x,s)\|_{L^{q}(I^h_l ,L^{r}(S))}ds.
\end{multline}
Applying Theorem \ref{thmstrichartz} for the linear semi-classical Schr\"odinger flow on $S$, the term to integrate in \eqref{vhlnorm} is bounded by $C \|\tilde{\Psi}(-h^{2}\Delta_{S})V_{h,l}(x,s)\|_{L^{2}(S)}$.
Using \cite[Thm.1.1]{ivpl08} and the fact that $\tilde{\chi}V_{h,l}=V_{h,l}$ (so that taking the norm over $\Omega$ or $S$  makes no difference) we obtain
\be\label{adem1}
\|v_{h,l,m}\|_{L^{q}(I^h_l ,L^{r}(\Omega))} \leq C\int_{I^h_l }\|V_{h,l}(x,s)\|_{L^{2}(\Omega)}ds.
\ee
After applying the Cauchy-Schwartz inequality in \eqref{estvhleeee}, \eqref{adem1} it remains to estimate the $L^{2}(I^h_l ,H^{\sigma}(\Omega))$ norm of $V_{h,l}$, where $\sigma\in\{0,n(\frac{1}{2}-\frac{1}{r})-\frac{1}{2}\}$. We do this using the precise form \eqref{defVhl1} and obtain 
\begin{multline}\label{estvhlmmee}
\|V_{h,l}\|_{L^{2}(I^h_l ,H^{\sigma}(\Omega))}\\
\leq C\|\varphi(t/h-l)[\Delta_{D},\chi]\Psi(-h^{2}\Delta_{D})e^{it\Delta_{D}}u_{0}\|_{L^{2}(I^h_l ,H^{\sigma}(\Omega))}\\
+Ch^{-1}\|\varphi'(t/h-l)\chi \Psi(-h^{2}\Delta_{D})e^{it\Delta_{D}}u_{0}\|_{L^{2}(I^h_l ,H^{\sigma}(\Omega))}.
\end{multline}
Since the operator $[\Delta_{D},\chi]\Psi(-h^{2}\Delta_{D})$ is bounded from $H^{\sigma+1}$ to $H^{\sigma}$, we deduce from \eqref{egalsd}, \eqref{estvhleeee}, \eqref{estvhlmmee}, \eqref{estvhlmm12} and Lemma \ref{lemschwartzrest} that\begin{multline}\label{estvhlmm12}
\|\tilde{\Psi}(-h^{2}\Delta_{D})v_{h,l}\|_{L^{q}(I^h_l ,L^{r}(\Omega))}\\ \leq Ch^{1/2}\|\tilde{\varphi}(t/h-l)\tilde{\chi}\Psi(-h^{2}\Delta_{D})e^{it\Delta_{D}}u_{0}\|_{L^{2}(I^h_l ,H^{1}_{0}(\Omega))}\\
+ Ch^{-1/2}\|\tilde{\varphi}(t/h-l)\tilde{\chi}\Psi(-h^{2}\Delta_{D})e^{it\Delta_{D}}u_{0}\|_{L^{2}(I^h_l ,L^{2}(\Omega))}\\
+C_{N}h^{N+1/2}\|\tilde{\varphi}(t/h-l)\tilde{\chi}\Psi(-h^{2}\Delta_{D})e^{it\Delta_{D}}u_{0}\|_{L^{2}(I^h_l ,H^{n(\frac{1}{2}-\frac{1}{r})+\frac{1}{2}}_{0}(\Omega))}\\
+C_{N}h^{N-1/2}\|\tilde{\varphi}(t/h-l)\tilde{\chi}\Psi(-h^{2}\Delta_{D})e^{it\Delta_{D}}u_{0}\|_{L^{2}(I^h_l ,H^{n(\frac{1}{2}-\frac{1}{r})-\frac{1}{2}}_{0}(\Omega))},
\end{multline}
where $\tilde{\varphi}\in C^{\infty}_{0}(\mathbb{R})$ is chosen equal to $1$ on the supports of $\varphi$.
Since $q\geq 2$ we estimate
\begin{multline}\label{estchi}
\|\tilde{\Psi}(-h^{2}\Delta_{D})\chi\Psi(-h^{2}\Delta_{D})e^{it\Delta_{D}}u_{0}\|^{q}_{L^{q}(\mathbb{R},L^{r}(\Omega))}\leq C \sum_{l=-\infty}^{\infty}\|\tilde{\Psi}(-h^{2}\Delta_{D})v_{h,l}\|^{q}_{L^{q}(I^h_l,L^{r}(\Omega))} \\
\leq Ch^{q/2}\Big(\sum_{l=-\infty}^{\infty}\|\tilde{\varphi}(t/h-l)\tilde{\chi}\Psi(-h^{2}\Delta_{D})e^{it\Delta_{D}}u_{0}\|^{2}_{L^{2}(I^h_l ,H^{1}_{0}(\Omega))}\Big)^{q/2}\\
{}+Ch^{-q/2}\Big(\sum_{l=-\infty}^{\infty}\|\tilde{\varphi}(t/h-l)\tilde{\chi}\Psi(-h^{2}\Delta_{D})e^{it\Delta_{D}}u_{0}\|^{2}_{L^{2}(I^h_l ,L^{2}(\Omega))}\Big)^{q/2}\\
{}+C_{N}h^{q(N+1/2)}\Big(\sum_{l=-\infty}^{\infty}\|\tilde{\varphi}(t/h-l)\tilde{\chi}\Psi(-h^{2}\Delta_{D})e^{it\Delta_{D}}u_{0}\|^{2}_{L^{2}(I^h_l ,H^{n(\frac{1}{2}-\frac{1}{r})+\frac{1}{2}}_{0}(\Omega))}\Big)^{q/2}\\
{}+C_{N}h^{q(N-1/2)}\Big(\sum_{l=-\infty}^{\infty}\|\tilde{\varphi}(t/h-l)\tilde{\chi}\Psi(-h^{2}\Delta_{D})e^{it\Delta_{D}}u_{0}\|^{2}_{L^{2}(I^h_l ,H^{n(\frac{1}{2}-\frac{1}{r})-\frac{1}{2}}_{0}(\Omega))}\Big)^{q/2}.
\end{multline}
The almost orthogonality of the supports of $\tilde{\varphi}(.-l)$ in time allows to estimate the term in the second line of \eqref{estchi} by
\be\label{boundbun}
Ch^{q/2}\|\tilde{\chi}\Psi(-h^{2}\Delta_{D})e^{it\Delta_{D}}u_{0}\|^{q}_{L^{2}(\mathbb{R},H^{1}_{0}(\Omega))},
\ee
the one in the third line by
\be\label{boundbun3}
Ch^{-q/2}\|\tilde{\chi}\Psi(-h^{2}\Delta_{D})e^{it\Delta_{D}}u_{0}\|^{q}_{L^{2}(\mathbb{R},L^{2}(\Omega))},
\ee
the term in the fourth line by
\be\label{boundbun4}
C_{N}h^{q(N+1/2)}\|\tilde{\chi}\Psi(-h^{2}\Delta_{D})e^{it\Delta_{D}}u_{0}\|^{q}_{L^{2}(\mathbb{R},H^{n(\frac{1}{2}-\frac{1}{r})+\frac{1}{2}}_{0}(\Omega))},
\ee
and the one in the last line of \eqref{estchi} by
\be\label{boundbunrest}
C_{N}h^{q(N-1/2)}\|\tilde{\chi}\Psi(-h^{2}\Delta_{D})e^{it\Delta_{D}}u_{0}\|^{q}_{L^{2}(\mathbb{R},H^{n(\frac{1}{2}-\frac{1}{r})-\frac{1}{2}}_{0}(\Omega))}.
\ee
We need the following smoothing effect on a non trapping domain:
\begin{prop}\label{propnic}( \cite[Prop.2.7]{bgt03})
Assume that $\Omega=\mathbb{R}^{n}\setminus\mathcal{O}$, where $\mathcal{O}\neq\emptyset$ is a compact non-trapping obstacle. Then for every $\tilde{\chi}\in C^{\infty}_{0}(\mathbb{R}^{n})$, $n\geq 2$, $\sigma\in [-1/2,1]$, one has
\be\label{smootheffect}
\|\tilde{\chi} \Psi(-h^{2}\Delta_{D})e^{it\Delta_{D}}u_{0}\|_{L^{2}(\mathbb{R},H^{\sigma+1/2}_{0}(\Omega))}\leq C\|\Psi(-h^{2}\Delta_{D})u_{0}\|_{H^{\sigma}(\Omega)}.
\ee
\end{prop}
\begin{rmq}
In \cite{bgt03}, Proposition \ref{propnic} is proved for $\sigma\in [0,1]$, but for spectrally localized data the result also follows using the estimates $(2.15)$ of \cite[Prop.2.7]{bgt03}.
\end{rmq}
We apply Proposition \ref{propnic} with $\sigma=1/2$ in \eqref{boundbun}, with $\sigma=-1/2$ in \eqref{boundbun3} and with $\sigma=n(\frac{1}{2}-\frac{1}{r})=\frac{2}{q}\in [0,1]$ in \eqref{boundbun4}. In \eqref{boundbunrest} we use that $n(\frac{1}{2}-\frac{1}{r})-\frac{1}{2}\leq \frac{1}{2}$ in order to estimate the $L^{2}(\mathbb{R}, H^{n(\frac{1}{2}-\frac{1}{r})-\frac{1}{2}}(\Omega))$ norm by the $L^{2}(\mathbb{R}, H^{\frac{1}{2}}(\Omega))$ norm and use Proposition \ref{propnic} with $\sigma=0$.  This yields
 \be
 \|\tilde{\Psi}(-h^{2}\Delta_{D})\chi\Psi(-h^{2}\Delta_{D})e^{it\Delta_{D}}u_{0}\|_{L^{q}(\mathbb{R},L^{r}(\Omega))}
 \leq C \|\Psi(-h^{2}\Delta_{D})u_{0}\|_{L^{2}(\Omega)},
 \ee
 where we used the spectral localization $\Psi$ to estimate $\|\Psi(-h^{2}\Delta_{D})u_{0}\|_{H^{\sigma}(\Omega)}$ by \\
 $h^{-\sigma}\|\Psi(-h^{2}\Delta_{D})u_{0}\|_{L^{2}(\Omega)}$.
This achieves the proof of Proposition \ref{propstri}.
 \end{itemize}
\end{proof}
In the rest of this section we show how Proposition \ref{propstri} implies Theorem \ref{thmstri}. We need the next lemma proved in \cite{ivpl08}:
\begin{lemma}\label{lemf}(see \cite[Thm.1.1]{ivpl08})
Let $\Psi_{0}\in C^{\infty}_{0}(\mathbb{R})$, $\Psi\in C^{\infty}_{0}((1/2,2))$ satisfy
\[
\Psi_{0}(\lambda)+\sum_{j\geq 1}\Psi(2^{-2j}\lambda)=1, \quad\forall\lambda\in \mathbb{R}.
\]
Then for all $r\in [2,\infty)$ we have
\be\label{besov}
\|f\|_{L^{r}(\Omega)}\leq C_{r}\Big(\|\Psi_{0}(-\Delta_{D})f\|_{L^{r}(\Omega)}+(\sum_{j=1}^{\infty}\|\Psi(-2^{-2j}\Delta_{D})f\|^{2}_{L^{r}(\Omega)})^{1/2}\Big).
\ee
\end{lemma}
Applying Lemma \ref{lemf} to $f=e^{it\Delta_{D}}u_{0}$ and taking the $L^{q}$ norm in time yields
\[
\|e^{it\Delta_{D}}u_{0}\|_{L^{q}(\mathbb{R},L^{r}(\Omega))}\leq\|\|e^{it\Delta_{D}}\Psi_{0}(-\Delta_{D})u_{0}\|_{L^{r}(\Omega)}+(\sum_{j\geq 1}\|e^{it\Delta_{D}}\Psi(-2^{-2j}\Delta_{D})u_{0}\|^{2}_{L^{r}(\Omega)})^{1/2}\|_{L^{q}(\mathbb{R})}
\]
which, by Minkowski inequality, leads to
$\|e^{it\Delta_{D}}u_{0}\|_{L^{q}(\mathbb{R},L^{r}(\Omega))}\leq C\|u_{0}\|_{L^{2}(\Omega)}$.
The proof of Theorem \ref{thmstri} is complete.

\subsection{Applications}
In this section we sketch the proofs of  Theorem \ref{thmgeq} and Theorem \ref{thmscattering}.

We start with Theorem \ref{thmgeq}. From Theorem \ref{thmstri} we have an estimate of the linear flow of the Schr\"{o}dinger equation
\be\label{lhs1}
\|e^{-it\Delta_{D}}u_{0}\|_{L^{5}(\mathbb{R},L^{30/11}(\Omega))}\leq C\|u_{0}\|_{L^{2}(\Omega)}.
\ee
One may shift regularity by $1$ and obtain
\be\label{lhs2}
\|e^{-it\Delta_{D}}u_{0}\|_{L^{5}(\mathbb{R},W^{1,30/11}(\Omega))}\leq C\|u_{0}\|_{H^{1}_{0}(\Omega)}.
\ee
Hence for small $T>0$ the left hand side in \eqref{lhs1}, \eqref{lhs2} will be small; for such $T$ let
$X_{T}:=L^{5}((0,T],W^{1,30/11}(\Omega))$. One may then set up the usual fixed point argument in $X_T$, as if $u\in X_{T}$ then $u^{5}\in L^{1}([0,T], H^{1}(\Omega))$.

Let us proceed with Theorem \ref{thmscattering}. From the work of Planchon and Vega \cite{plve08}, one has a global in time control  on the solution $u$, at the level of  $\dot H^\frac 1 4$ regularity: 
$$
u\in L^4((0,+\infty),L^4(\Omega)).
$$
By interpolation with either mass or energy conservation, combined with the local existence theory, one may bootstrap this global in time control into
$$
u\in L^{p-1}((0,+\infty),L^\infty(\Omega)),
$$
from which scattering in $H^1_0(\Omega)$ follows immediately.

\section{Appendix}
\subsection{Finite speed of propagation for the semi-classical equation}
In this section we recall several properties of the semi-classical Schr\"odinger flow (for further discussions and proofs we refer the reader to \cite{le92}). Let $S$ be a compact manifold with smooth boundary $\partial S$.

\begin{dfn}
We say that a symbol $q(y,\eta)\in S^{m}_{\rho,\delta}$ is of type $(\rho,\delta)$ and of order $m$ if 
\[
\forall\alpha,\beta \quad \exists C_{\alpha,\beta}>0\quad |\partial^{\beta}_{y}\partial^{\alpha}_{\eta}q(y,\eta)|\leq C_{\alpha,\beta}(1+|\eta|)^{m-\rho|\alpha|+\delta|\beta|}.
\]
\end{dfn}
For $q\in S^{m}_{1,0}$ we let $Op_{h}(q)=Q(y,hD,h)$ be the $h$-pseudo differential operator defined by
\[
Op_{h}(q)f(y)=\frac{1}{(2\pi h)^{n}}\int e^{\frac{i}{h}(y-\tilde{y})\eta}q(y,\eta,h)f(\tilde{y})d\tilde{y}.
\]
We set $y=(x,t)\in S\times\mathbb{R}$ and denote $\eta=(\xi,\tau)$ the dual variable of $y$. Near a point $x_{0}\in\partial S$ we can choose a system of local coordinates such that $S$ is given by $S=\{x=(x',x_{n})|x_{n}>0\}$. We define the tangential operators 
\[
Op_{h,\text{tang}}(q)f(y)=\frac{1}{(2\pi h)^{n-1}}\int e^{\frac{i}{h}(y'-\tilde{y}')\eta'}q(y,\eta',h)f(\tilde{x}',x_{n},\tilde{t})d\tilde{y}'d\eta',
\] 
where $y=(x',x_{n},t)$, $y'=(x',t)$, $\tilde{y}'=(\tilde{x}',\tilde{t})$, $\eta=(\xi',\xi_{n},\tau)$, $\eta'=(\xi',\tau)$ and where the symbol $q(y,\eta',h)\in S^{m}_{1,0,\text{tang}}$ i.e. such that
\[
\forall\alpha,\beta \quad \exists C_{\alpha,\beta}>0\quad |\partial^{\alpha}_{y}\partial^{\beta}_{\eta'}q(y,\eta',h)|\leq C_{\alpha,\beta}(1+|\eta'|)^{m-|\beta|}.
\]
\bigskip

In what follows we let $(S,g)$ be a compact Riemannian manifold with strictly concave boundary satisfying the Assumptions \ref{assum1}. Let also $v_{0}\in L^{2}(S)$ be compactly supported outside a small neighborhood of the boundary, $\Psi\in C^{\infty}_{0}((\alpha_{0},\beta_{0}))$ and let $v(x,t)=e^{iht\Delta_{g}}\Psi(-h^{2}\Delta_{g})v_{0}$ denote the linear semi-classical Schr\"odinger flow with initial data at time $t=0$ equal to $\Psi(-h^{2}\Delta_{g})v_{0}$ and such that $\|\Psi(-h^{2}\Delta_{g})v_{0}\|_{L^{2}(S)}\lesssim 1$.
 
Let $\pi:T^{*}(\bar{S}\times\mathbb{R})\rightarrow T^{*}(\partial S\times\mathbb{R})\cup T^{*}(S\times\mathbb{R})$ be the canonical projection defined, for $y=(x,t)$, $\eta=(\xi,\tau)$ by
\[
\pi|_{T^{*}(S\times\mathbb{R})}=Id,\quad \pi(y,\eta)=(y,\eta|_{T^{*}(\partial S\times\mathbb{R})}),\quad \ for \ (y,\eta)\in T^{*}(\bar{S}\times\mathbb{R})|_{\partial S\times\mathbb{R}}.
\]
We introduce the characteristic set 
\[
\Sigma_{b}:=\pi\{(y,\eta)|\eta=(\xi,\tau), \tau+|\xi|^{2}_{g}=0, -\beta_{0}\leq\tau\leq -\alpha_{0}\},
\]
where $|\xi|^{2}_{g}=<\xi,\xi>_{g}=:\xi^{2}_{n}+r(x,\xi')$ denotes the inner product given by the metric $g$ and where, due to the strict concavity of the boundary we have $\partial_{x_{n}}r(x',0,\eta')<0$.
\begin{dfn}\label{dfnwf}
We say that a point $\rho_{0}=(y_{0},\eta_{0})\in T^{*}_{b}(\partial S\times\mathbb{R}):=T^{*}(\partial S\times\mathbb{R})\cup T^{*}(S\times\mathbb{R})$ doesn't belong to the $b$-wave front set $WF_{b}(v)$ of $v$ if there exists a $h$-pseudo-differential operator of symbol $q(y,\eta,h)$ (respectively $q(y,\eta',h)$ if $\rho_{0}\in T^{*}(\partial S\times\mathbb{R})$) with compact support in $(y,\eta)$, elliptic at $\rho_{0}$, and a smooth function $\phi\in C^{\infty}_{0}$ equal to $1$ near $y_{0}$, such that for every $\sigma\geq 0$ the following holds
\[
\forall N\geq 0\quad \exists C_{N}>0\quad \|Op_{h}(q)\phi v\|_{H^{\sigma}(S\times\mathbb{R})}\leq C_{N}h^{N}.
\]
We shall write $\rho_{0}\notin WF_{b}(v)$.
\end{dfn}
\begin{prop}\label{propellreg}(Elliptic regularity \cite[Thm.3.1]{le92})
Let $q(y,\eta)$ a symbol such that $q=0$ on a neighborhood of $\Sigma_{b}$. Then for every $\sigma\geq 0$ we have
\[
 \forall N\geq 0\quad \exists C_{N}>0\quad \|Op_{h}(q)v\|_{H^{\sigma}(S)}\leq C_{N}h^{N}.
\]
\end{prop}
Proposition \ref{propellreg} is proved by Lebeau \cite{le92} for eigenfunctions of the Laplace operator, but the same arguments apply in this setting.
From Proposition \ref{propellreg} and \cite[Sections 2,3]{le92} it follows:
\begin{cor}\label{corwf}
There exists a constant $D>0$ such that
\[
WF_{b}(v)\subset\Sigma_{b}\cap \{-\tau\in [\alpha_{0},\beta_{0}], |\xi|_{g}\leq D\}.
\]
\end{cor}
\begin{cor}\label{cort}(\cite[Chp.3]{le92})
Let $\varphi\in C^{\infty}_{0}(\mathbb{R})$ be equal to $1$ near the interval $[-\beta_{0},-\alpha_{0}]$. Then for $t$ in any bounded interval $I$ one has
\be\label{llochdt}
\forall N\geq 1,\quad \exists C_{N}>0\quad |(1-\varphi)(hD_{t})v|\leq C_{N}h^{N}, \quad \forall t\in I.
\ee
\end{cor}
\begin{cor}\label{corregel}(Elliptic regularity at "$\infty$")
Let $\vartheta\in C^{\infty}_{0}(\mathbb{R}^{n})$ be equal to $1$ on $\{|\xi|_{g}\leq D\}$, then
\be\label{lelochdt}
\forall N\geq 1,\quad \exists C_{N}>0\quad |(1-\vartheta)(hD_{x})v|\leq C_{N}h^{N}.
\ee
\end{cor}
\begin{prop}\label{prople92}
Let $y_{0}\notin pr_{y}(WF_{b}(v))$, where by $pr_{y}$ we mean the projection on the variable $y=(x,t)$. Then there exists $\phi\in C^{\infty}_{0}$, $\phi=1$ near $y_{0}$ such that for every $\sigma\geq 0$ we have
\[
\forall N\geq 0 \quad \exists C_{N}>0\quad \|\phi v\|_{H^{\sigma}(S)}\leq C_{N}h^{N}.
\]
\end{prop}
\begin{proof}
Let $\varphi$, $\vartheta$ be the functions defined in Corollaries \ref{cort}, \ref{corregel}. Then, using again Proposition \ref{propellreg} we have
\be\label{suppdualv}
v(x,t)=\varphi(hD_{t})\vartheta(hD_{x})v+O(h^{\infty}).
\ee
Let now $y_{0}=(x_{0},t_{0})\notin pr_{y}(WF_{b}(v))$. It follows that for every $\eta\neq 0$, $(y_{0},\eta)\notin WF_{b}(v)$ and in particular for every $\eta_{0}\in \text{supp}(\vartheta)\times\text{supp}(\varphi)$ there exists a symbols $q_{0}(y,\eta,h)$ with compact support in $(y,\eta)$ near $(y_{0},\eta_{0})$ and elliptic at $(y_{0},\eta_{0})$, and there exists $\phi_{0}\in C^{\infty}_{0}$ equal to $1$ in a neighborhood $U_{0}$ of $y_{0}$ such that for every $\sigma\geq 0$
\[
\forall N\geq 0\quad \exists C_{N}>0\quad \|Op_{h}(q_{0})\phi v\|_{H^{\sigma}(S)}\leq C_{N}h^{N}.
\]
Eventually shrinking $U_{0}$, suppose that $q_{0}$ is elliptic on $U_{0}\times W_{0}$ where $W_{0}$ is an open neighborhood of $\eta_{0}$. Then it follows that on $U_{0}$, for every $\sigma\geq 0$
\[
\forall N\geq 0\quad \exists C_{N}>0\quad \|\phi v\|_{H^{\sigma}(U_{0})}\leq C_{N}h^{N}.
\]
Since the set $\text{supp}(\vartheta)\times\text{supp}(\varphi)$ is compact there exist $\eta^{\alpha}$, $\alpha\in\{1,..,N\}$ for some fixed $N\geq 1$ and for each $\eta^{\alpha}$ there exist symbols $q_{\alpha}$ elliptic on some neighborhoods $U_{\alpha}\times W_{\alpha}$ of $(y_{0},\eta^{\alpha})$ and smooth functions $\phi_{\alpha}\in C^{\infty}_{0}$ equal to $1$ on the neighborhoods $U_{\alpha}$ of $y_{0}$, such that $\text{supp}(\vartheta)\times\text{supp}(\varphi)\subset \cup_{j=1}^{N}W_{\alpha}$. Let $\phi\in C^{\infty}_{0}$ be equal to $1$ in an open neighborhood of $y_{0}$ strictly included in the intersection $\cap_{\alpha=1}^{N}U_{\alpha}$ (which has nonempty interior) and supported in the compact set $\cap_{\alpha=1}^{N}\text{supp}(\phi_{\alpha})$. Considering a partition of unity associated to $(U_{\alpha}\times W_{\alpha})_{\alpha}$ and using \eqref{suppdualv} we deduce that $\phi$ satisfies Proposition \ref{prople92}. 
\end{proof}
\begin{prop}\label{propschtdd}\cite[Lemma B.7]{NBth}
Let $v(x,t)=e^{ith\Delta_{g}}\Psi(-h^{2}\Delta_{g})v_{0}$ like before, $v_{0}\in L^{2}(S)$ and let $Q$ be a $h$-pseudo-differential operator of order $0$, $t_{0}>0$ and $\tilde{\psi}\in C^{\infty}_{0}((-2t_{0},-t_{0}))$. Let $w$ denote the solution to
\be
\left\{\begin{array}{ll}
(ih\partial_{t}+h^{2}\Delta_{g})w=ih\tilde{\psi}(t)Q(v),\quad \text{on}\quad S\times\mathbb{R},\\
w|_{\partial S}=0,\quad w|_{t<-2t_{0}}=0.
\end{array}
\right.
\ee
If $\rho_{0}\in WF_{b}(w)$ then the broken bicharacteristic starting from $\rho_{0}$ has a nonempty intersection with $WF_{b}(v)\cap \{t\in\text{supp}(\tilde{\psi})\}$.
\end{prop}

\subsection{Proof of Lemma \ref{lemschwartzrest}}
In this section $(M, \Delta_{M})$ denotes either $(S,\Delta_{S})$ or $(\Omega,\Delta_{D})$, respectively. This notation will be used to refer both domains at the same time. Let $\tilde{\chi}\in C^{\infty}_{0}(\mathbb{R}^{n})$ be such that $\Delta_{D}\tilde{\chi}=\Delta_{S}\tilde{\chi}$.

Let $\varphi_{0}\in C^{\infty}(\mathbb{R})$ be supported in the interval $[-4,4]$ and $\varphi\in C^{\infty}(\mathbb{R})$ be supported in $[-4,-1]\cup [1,4]$ such that for all $\xi\in\mathbb{R}$
\[
\varphi_{0}(\xi)+\sum_{k\geq 1}\varphi(2^{-k}\xi)=1.
\]
If $\hat{\Psi}$ denotes the Fourier transform of $\Psi$, we write it using the preceding sum 
\[
\hat{\Psi}(\xi)=\hat{\Psi}(\xi)\Big(\varphi_{0}(\xi)+\sum_{k\geq 1}\varphi(2^{-k}\xi)\Big)
\]
and denote $\phi_{k}\in \mathcal{S}(\mathbb{R})$ the functions such that $\hat{\phi}_{0}(\xi)=\hat{\Psi}(\xi)\varphi_{0}(\xi)$, $\hat{\phi}_{k}(\xi)=\hat{\Psi}(\xi)\varphi(2^{-k}\xi)$. We denoted by $\mathcal{S}(\mathbb{R})$ the Schwartz space of rapidly decreasing functions. Hence we have
\be\label{decompsi}
\Psi(\lambda)=\sum_{k\in\mathbb{N}}\phi_{k}(\lambda),\quad\text{where}\quad \|\hat{\phi}_{k}\|_{L^{\infty}}=\|\hat{\Psi}(\xi)\varphi(2^{-k}\xi)\|_{L^{\infty}}\leq C_{N}2^{-kN},\quad \forall N\in \mathbb{N}.
\ee
For $k\in\mathbb{N}$ write
\be\label{phikf1}
\phi_{k}(h\sqrt{-\Delta_{M}})\tilde{\chi}v_{h,l}=\frac{1}{2\pi}\int_{\text{supp}(\hat{\phi}_{k})} e^{i\xi h\sqrt{-\Delta_{M}}}\tilde{\chi}v_{h,l}\hat{\phi}_{k}(\xi)d\xi.
\ee
On the support of $\hat{\phi}_{k}(\xi)$, $|\xi|\simeq 2^{k}$ and for $k\leq \frac{1}{2}\log_{2}(1/h)$ for example we see, by the finite speed of propagation of the wave operator, that
on a time interval of size $2^{k}h \leq h^{1/2}$ we remain in a fixed neighborhood
of the boundary of $\Omega$ where $\Delta_{D}$ coincides with $\Delta_{S}$, therefore we can introduce
$\chi_{1}$ equal to $1$ on a fixed neighborhood of the support of $\tilde{\chi}$ (independent of $k$, $h$) such that for every $k\leq \frac{1}{2}\log_{2}(1/h)$
\begin{equation}\label{restschw}
\chi_{1}\phi_{k}(h\sqrt{-\Delta_{S}})\tilde{\chi}v_{h,l}
=\chi_{1}\phi_{k}(h\sqrt{-\Delta_{\Omega}})\tilde{\chi}v_{h,l}.
\end{equation}
Since $v_{h,l,s}=\chi_{1}(\tilde{\Psi}(-h^{2}\Delta_{D})-\tilde{\Psi}(-h^{2}\Delta_{S}))v_{h,l}$ and $v_{h,l}=\tilde{\chi}v_{h,l}$, we obtain, using \eqref{restschw}
\be
v_{h,l,s}=\chi_{1}\Big(\sum_{k\geq \frac{1}{4}\log_{2}(1/h)}(\phi_{k}(h\sqrt{-\Delta_{\Omega}})-\phi_{k}(h\sqrt{-\Delta_{S}})) \Big)\tilde{\chi}v_{h,l}.
\ee
In order to estimate the $L^{q}(I^h_l ,L^{r}(\Omega))$ norm of $v_{h,l,s}$ it will be enough to estimate separately the norms of  $\chi_{1}\phi_{k}(h\sqrt{-\Delta_{M}})\tilde{\chi}v_{h,l}$ for $k\geq \frac{1}{4}\log_{2}(1/h)$ where $(M,\Delta_{M})\in \{(\Omega,\Delta_{D}), (S,\Delta_{S})\}$. Using the Cauchy-Schwartz inequality and the Sobolev embeddings gives
\begin{multline}\label{hohoest}
\|\chi_{1}\phi_{k}(h\sqrt{-\Delta_{M}})\tilde{\chi}v_{h,l}\|_{L^{q}(I^h_l ,L^{r}(\Omega))}\leq Ch^{1/q}\|\chi_{1}\phi_{k}(h\sqrt{-\Delta_{M}})\tilde{\chi}v_{h,l}\|_{L^{\infty}(I^h_l ,L^{r}(\Omega))}\\
\leq Ch^{1/q}\|\chi_{1}\phi_{k}(h\sqrt{-\Delta_{M}})\tilde{\chi}v_{h,l}\|_{L^{\infty}(I^h_l ,H^{n(\frac{1}{2}-\frac{1}{r})}(\Omega))}\\
\leq C_{N}h^{1/q}2^{-kN}\|\tilde{\chi}v_{h,l}\|_{L^{\infty}(I^h_l ,H^{n(\frac{1}{2}-\frac{1}{r})}(\Omega))},\quad \forall N\in \mathbb{N},
\end{multline}
where in the last line we used \eqref{decompsi}.
We estimate the last term in \eqref{hohoest} writing the Duhamel formula for $v_{h,l}$ only on $\Omega$ using the equation \eqref{eqvh}, since in this case the smoothing effect yields (see \cite{stta02}, \cite{bgt03} or the dual estimates of \eqref{smootheffect} in Proposition \ref{propnic}):
\be
\|\tilde{\chi}v_{h,l}\|_{L^{\infty}(I^h_l ,H^{n(\frac{1}{2}-\frac{1}{r})}(\Omega))}\leq C\|V_{h,l}\|_{L^{2}(I^h_l ,H^{n(\frac{1}{2}-\frac{1}{r})-\frac{1}{2}}(\Omega))}.
\ee
Since we consider here only large values $k\geq \frac{1}{4}\log_{2}(1/h)$, each $2^{-k}$ is bounded by $h^{1/4}$, therefore, after summing over $k$ we obtain
\be
\|v_{h,l,s}\|_{L^{q}(I^h_l ,L^{r}(\Omega))}\leq C_{N}h^{1/q+N/4}\|V_{h,l}\|_{L^{2}(I^h_l ,H^{n(\frac{1}{2}-\frac{1}{r})-\frac{1}{2}}(\Omega))},\quad \forall N\in \mathbb{N}.
\ee

\bibliography{secbib}

\end{document}